\documentclass[11pt,reqno,oneside]{amsart}

\usepackage{amsmath} 
\usepackage{amssymb}
\usepackage{graphicx} \DeclareGraphicsExtensions{.eps} 
\usepackage{caption}
\usepackage{placeins}
\usepackage{amsbsy} 
\usepackage[utf8]{inputenc}
\usepackage{bm}

\newcommand{\dx}{\textrm{d}x} 
\newcommand{\dy}{\textrm{d}y}

 \numberwithin{equation}{section}

\newtheorem{thm}{Theorem}[section] 
\newtheorem{prop}[thm]{Proposition}
\newtheorem{lemma}[thm]{Lemma}

\newtheorem{defi}{Definition}[section]

\newtheorem{rmk}[thm]{Remark}

\newcommand{\R}{\mathbb{R}} 
\newcommand{\N}{\mathbb{N}}

\newcommand{\diam}{\textrm{diam}} 

\newcommand{\dt}{\textrm{d}t} 
\newcommand{\ds}{\textrm{d}s}
\newcommand{\dr}{\textrm{d}r} 
\newcommand{\dtheta}{\textrm{d}\theta}
\newcommand{\dxi}{\textrm{d}\xi}

\begin{document}

\author{Ignacio Ojea} \address{Departamento de Matem\'atica Facultad de Ciencias
Exactas y Naturales, Universidad de Buenos Aires and IMAS - CONICET, 1428,
Buenos Aires, Argentina} \email{iojea@dm.uba.ar} \title{Anisotropic
regularity for elliptic problems with Dirac measures as data}

\begin{abstract} 
We study the Possion problem with singular data given by a source supported on a one dimensional curve strictly contained in a three dimensional domain. We prove regularity results for the solution on isotropic and on anisotropic weighted spaces of Kondratiev type. Our technique is based on the study of a regularized problem. This allows us to exploit the local nature of the singularity. Our results hold with very few smoothness hypotheses on the domain and on the support of the data. We also discuss some extensions of our main results, including the two dimensional case, sources supported on closed curves and on polygonals. 
\end{abstract}

\maketitle

\noindent{\bf Keywords:} Anisotropy, Dirac delta, singular data, Weighted Sobolev
spaces

\noindent{\bf Funding:} This work was supported by ANPCyT under grant PICT 2018 - 3017, by CONICET under grant PIP112201130100184CO and by Universidad de Buenos Aires under grant 20020170100056BA. 
\section{Introduction} 
In this paper we study the regularity in isotropic and anisotropic weighted spaces
of the solution of the problem: 
\begin{equation}\label{problem} \left\{
		\begin{array}{cl} -\Delta u = \sigma\delta_\Lambda &  \textrm{ in }\Omega,\\ u = 0
& \textrm{ in }\partial\Omega, \end{array}\right. \end{equation} 
where $\Omega$ is a bounded domain in $\R^3$, $\Lambda$ is a curve strictly
contained in $\Omega$ and $\sigma$ is a function defined over $\Lambda$.
$\delta_\Lambda$ is a Dirac delta supported on $\Lambda$. 

Such problems arise in fluid mechanics, as a simplificaction of a complex system used for saving computational resources, see for example \cite{DQ_1d3d,D_1d3d}. The singular data can also be used for modelling an idealized load supported on $\Lambda$. 

Our results are inspired by \cite{DQ_1d3d,D_1d3d} where the coupling of two
diffusion-reaction problems (one in $3D$, the other in $1D$) is studied as a
model for blood flow through tissue. There, a certain regularity of the
solution in weighted Sobolev spaces is assumed in order to obtain error
estimates for the approximation of $u$ via a finite element method. Such
regularity was later proven in \cite{AdCN1,AdCN2}. Our goal is to extend the results
of these articles taking into account the anisotropic behaviour of the solution. 

In \cite{AdCN1} only the case where $\Lambda$ is a straigth line is considered, and a technique based on Fourier and Mellin transforms is applied. In \cite{AdCN2} the singularity is supported on a curve, but a smooth transformation is applied in order to straighten it. Then, the isotropic weighted regularity is obtained by a technique based on a priori estimates proven on a dihedron or a cone with singularities. These estimates are obtained for differential operators with variable coefficients, which arise as a consequence of the straightening of $\Lambda$. In this context, two main assumptions are imposed in order to obtain regularity results for derivatives of $u$ of order $m$: $\Lambda$ is assumed to be of class $\mathcal{C}^{m+2}$ and $\Omega$ is assumed to be of class $\mathcal{C}^{m}$.

Our approach is based on the regularization of the data. For each $\rho>0$ we define a smooth function $\sigma_\rho$ supported on a neighbourhood of $\Lambda$ and such that $\sigma_\rho\to \sigma\delta_\Lambda$ in a distributional sense when $\rho\to0$. We then study the solution $u_\rho$ of the regularized problem $-\Delta u_\rho=\sigma_\rho$. In particular, we consider weighted norms of $u_\rho$ and its derivates and establish conditions on the weights that allows us to take limit with $\rho$ tending to $0$, thus obtaining regularity results for the solution $u$ of the singular problem \eqref{problem}. This method is local in nature and has some advantages with respect to the ideas applied in \cite{AdCN1,AdCN2}. On the one hand we only need $\Lambda$ to be smooth enough so that curvilinear cylindrical coordinates can be defined in a neighbourhood of it. On the other hand, since the singularity of the data is localized at $\Lambda$, which is far from $\partial\Omega$, we only need the domain to be regular enough so that it does not introduce new singularities. In general, in order to obtain estimates for the derivatives up to order $m$ we assume that $\Omega$ is of class $\mathcal{C}^{m-1,1}$. However, this can be relaxed in some particular cases. For example: it is well known that if $\Omega$ is a convex polyhedron the solution of the Laplace equation with regular data belongs to $H^2(\Omega)$. We take advantage of this fact for proving that our results stands for derivatives of order $m\le 2$ on convex polyedra.  

Furthermore the same regularization technique can be applied for obtaining regularity results on \emph{anisotropic} weighted spaces, assuming regularity of the derivatives of $\sigma$ along $\Lambda$. In this case we assume for convenience that $\Lambda$ is a straight segment. Near the center of $\Lambda$ the derivatives of $u$ in direction parallel to $\Lambda$ are smoother than the ones in other directions. However, a sharper singularity arises at the extreme points of the curve, so the regularity is proven in weighted spaces involving two weights: one given by a power of the distance to $\Lambda$ and another given by a power of the distance to its extreme points.  Our anisotropic result resembles well known regularity results for elliptic problems on polyhedral domains with interior edges, where singularities arise at the interior edges and at the vertices adyacent to them. 

Another interesting feature of our approach is that it can be applied, with little adaptations, to some special cases. In Section \ref{section extensions} we discuss some of them. In particular: even though we treat extensively the case where $\Lambda$ is an open simple curve, it is easy to extend our results to \emph{closed} simple curves. Moreover, a version of our anisotropic results can be obtained when $\Lambda$ is a polygonal fracture. In that case, the vertices of the polygonal $\Lambda$ act as extreme points of the segments that form $\Lambda$. Finally, we also comment the two dimensional case, where the same ideas can be applied.


\section{Preliminaries}\label{section preliminaries}

In the sequel, $C$ denotes a constant that may change from line to line. When relevant, we indicate the dependance of $C$. For example: $C(\gamma)$ is a constant depending on the parameter $\gamma$.  We say that two quantities $a$ and $b$ are proportional, and we denote $a\sim b$ if there are constants $C_1$ and $C_2$ such that $C_1 a\le b\le C_2 a$. For every set $E$ we denote $|E|$ the measure of $E$, $|E|=\int_E\dx$. Moreover, $\chi_E(x)$ stands for the characteristic function of $E$, which takes the value $1$ for $x\in E$ and vanishes outside $E$. Given an exponent $1<p<\infty$, $p'$ stands for its H\"older conjugate: $1/p+1/p'=1$.

We denote $\beta=(\beta_1,\beta_2,\beta_3)\in \mathbb{N}_0^3$ a multiindex and $|\beta|=\beta_1+\beta_2+\beta_3$, its order. $D^\beta u$ stands for the derivative $\partial^{\beta_3}_{x_3}\partial^{\beta_2}_{x_2}\partial^{\beta_1}_{x_1} u$. 

We consider $\Lambda\subset \Omega$ a simple curve given by: \[\Lambda
=\{x\in\R^3:\,x=\lambda(s): s\in[0,L]\},\] where $s$ is a curvilinear abscissa and $\lambda$ is a smooth parametrization by arc-length. For each $s\in[0,L]$, we denote ${\bf t}(s)$, ${\bf n}(s)$ and ${\bf b}(s)$ the tangent, normal an binormal versors on $\Lambda$. For every $\rho>0$ we consider a cylindrical neighbourhood of $\Lambda$ given by: 
\[C(\Lambda,\rho) = \{x\in\R^3:\,
x=X(r,\theta,s), \, (r,\theta,s)\in [0,\rho)\times[0,2\pi)\times[0,L]\}\] where:
\[X(r,\theta,s) = \lambda(s) + r \cos(\theta){\bf n}(s) + r \sin(\theta){\bf
b}(s).\]

The case of a closed curve, where $\lambda(0)=\lambda(L)$ is briefly considered in Section \ref{section extensions}. For now, let us assume that $\Lambda$ is not closed. In this case, we will also need to consider neighbourhoods of the endpoints. We define:
\begin{align*}
	B^+_\rho(0) &= \{\lambda(0)+r\cos(\theta)\sin(\xi){\bf n}(0)+r\sin(\theta)\sin(\xi){\bf b}(0),\; (r,\theta,\xi)\in [0,\rho)\times[0,2\pi)\times(\pi,2\pi)\},\\
	B^+_\rho(L) &= \{\lambda(L)+r\cos(\theta)\sin(\xi){\bf n}(0)+r\sin(\theta)\sin(\xi){\bf b}(0),\; (r,\theta,\xi)\in [0,\rho)\times[0,2\pi)\times(0,\pi)\}.
\end{align*}
$B^+_\rho(0)$ and $B^+_\rho(L)$ are half spheres around $\lambda(0)$ and $\lambda(L)$ respectively, but outside $C(\Lambda,\rho)$. 

Finally, let us denote: 
\[B(\Lambda,\rho) = \{x\in\R^3:\; d(x,\Lambda)<\rho\}.\]

\begin{figure}[!h]
	\includegraphics[scale=0.6]{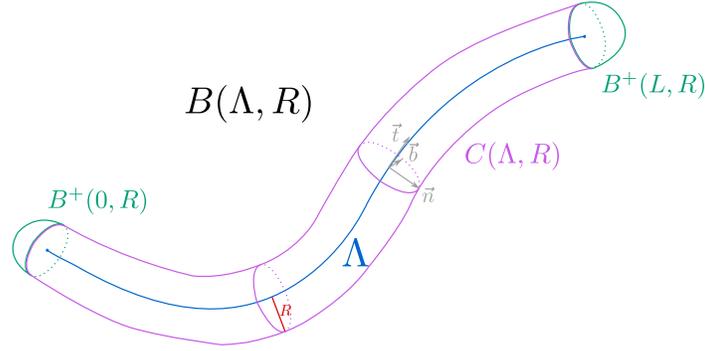}
	\caption{A curve $\Lambda$ and a neighbourhood $B(\Lambda,R)$.}
	\label{fig ball around curve}
\end{figure}

Then, we have that if $\Lambda$ is an open curve: 
\[B(\Lambda,\rho) = B^+_\rho(0) \cup C(\Lambda,\rho)\cup B^+_\rho(L). \]
If $\Lambda$ is closed, $B^+_\rho(0)$ and $B^+_\rho(L)$ are unnecessary and $B(\Lambda,\rho)=C(\Lambda,\rho)$. Figure \ref{fig ball around curve} shows an example of $B(\Lambda,R)$ for an open curve.

We assume that $\Lambda$ is smooth enough so there is a radius $R_0$ such that
$B(\Lambda,R_0)\subset \Omega$ and the projection from $B(\Lambda,R_0)$ to
$\Lambda$ is unique, i.e.: \[\forall x\in B(\Lambda,R_0)\,\exists!
s_x\in[0,L]:\; d(x,\Lambda)=\|x-\lambda(s_x)\|.\]
For $x\in B^+(0,R_0)$ and $x\in B^+(L,R_0)$, the distance is reached at an extreme point of $\Lambda$, i.e.: $s_x=0$ and $s_x=L$ respectively. On the other hand, for $x\in C(\Lambda,R_0)$, we have that $|x-\lambda(s_x)|=r_x$, the radial component of the cylindrical coordinates defined by $X$.  

We also assume that $\sigma\in L^2(\Lambda)$. Moreover, in order to simplify the notation, we identify $\Lambda$ with the interval $[0,L]$ and we write $\sigma(t)$ instead of $\sigma(\lambda(t))$. 

Since the solution $u$ does not belong to $H^1$, we need to study problem
\eqref{problem} in a non-standard setting. We work in weighted Sobolev and
Kondratiev-type spaces. Given $\omega$ a non-negative function defined on
$\Omega$, we denote $L^p(\Omega,\omega)$ the space of functions $v$ such that
$\|v\omega^\frac{1}{p}\|_{L^p(\Omega)}<\infty$. Our results are stated for $p=2$, but
other values of $p$ are considered in some technical arguments.
$H^k(\Omega,\omega)$ is the space of functions in $L^2(\Omega,\omega)$ with
weak derivatives up to order $k$ in $L^2(\Omega,\omega)$, and
$H^k_{0}(\Omega,\omega)$ is the closure of $C^\infty_0(\omega)$ in
$H^k(\Omega,\omega)$. 

We denote $d(x)=d(x,\Lambda)$, the distance from $x$ to $\Lambda$. Our isotropic results are given for weights of
the form: $d(x)^\gamma$, so we simplify the notation defining
$L^p_\gamma(\Omega)=L^p(\Omega,d^{p \gamma})$, $H^1_\gamma(\Omega)=
H^1(\Omega,d^{2\gamma})$ and $H^1_{0,\gamma}(\Omega)=
H^1_0(\Omega,d^{2\gamma})$. It is important to notice that
$L^2_{\gamma}(\Omega)\subset L^2_{\gamma+\mu}(\Omega)$ with continuity
for every $\mu>0$. 



As a consequence of \cite[Lemma 3.3]{DLG_Stokes} (see also \cite{Dyda19}) we have that if
$-1<\gamma<1$, $d(x)^{2\gamma}$ belongs to the Muckenhoupt class $A_2$. This
implies that the Rellich-Kondrakov theorem and the Poincaré inequality hold
on $H^1_{\gamma}$.

Our first goal is to give, for some values of $\gamma$ a weighted setting
for problem \eqref{problem} of the form:
\begin{equation}\label{weakproblem}
\textrm{Find }u\in H^1_{0,\gamma}(\Omega):\quad \int_\Omega \nabla u \nabla
v = \int_\Lambda \sigma v,\quad \forall v\in H^1_{0,-\gamma}(\Omega). 
\end{equation}

The first step is to prove that the right-hand side is well defined. In
\cite{DQ_1d3d} it is proven that for $0<\gamma<1$ there is a unique continuous
trace operator $T_\Lambda:H^1_{-\gamma}(\Omega)\to L^2(\Lambda)$. Here we
apply essentially the same argument for proving that the measure
$\sigma\delta_\Lambda$ is a bounded operator on $H^1_{-\gamma}(\Omega)$. For this,
we need the following weighted Hardy inequality (see \cite[page 6]{KPS_Hardy}
and \cite[Section 1]{KO_Hardy}.):

\begin{thm}[Weighted Hardy inequality] Let $0<p\le q<\infty$, $0<r\le\infty$ and $\omega$ and $\nu$ be weight functions defined on $(0, \infty)$. Assume that, for every $r > 0$, 
	\[\int_0^r \nu(t)^\frac{1}{1-p}\dt <\infty.\] 
Then, the inequality 
\begin{equation}\label{Hardy} 
	\left(\int_0^R \bigg(\int_0^r f(t) \dt \bigg)^q \omega(r) \dr \right)^\frac{1}{q} \le C \left(\int_0^R f(t)^p \nu(t) \dt \right)^\frac{1}{p}, 
\end{equation} 
holds for every positive function $f$ on $(0,\infty)$ if and only if:  
\[D = \sup_{r\in(0,R)} \left(\int_r^R\omega(t)\dt\right)^\frac{1}{q} 
\left(\int_0^r\nu(t)^\frac{1}{1-p}\dt\right)^\frac{p-1}{p}<\infty.\]
Moreover, the best constant in \eqref{Hardy} satisfies the estimate 
\[D \le C \le k(p, q)D,\] where  
\[k(p, q) = \left(\frac{p+qp-q}{p}\right)^\frac{1}{q}
\left(\frac{p+qp-q}{(p-1)q}\right)^\frac{p-1}{p}.\]
\end{thm}

\begin{thm}\label{theorem continuity delta} 
	If $\sigma\in L^2(\Lambda)$ and $0<\gamma<1$, we have that $\sigma\delta_\Lambda \in (H^1_{-\gamma}(\Omega))'$, and the following estimate holds: 
	\[|(\sigma\delta_\Lambda)(v)| \le \sqrt{\frac{2}{\pi}} \|\sigma\|_{L^2(\Lambda)}\left(R^{-2}\|v\|^2_{L^2(C(\Lambda,R))} + R^{2\gamma}c(\gamma)\|\nabla v\|^2_{L^2_{-\gamma}(C(\Lambda,R))} \right)^\frac{1}{2},\] 
where $c(\gamma)$ is a constant that tends to $\infty$ as $\gamma\to 0$. \end{thm} 
\begin{proof} 
	By a density argument, it is enough to prove the result for every $v\in C^\infty(\Omega)$. We have
that $(\sigma\delta_\Lambda)(v) = \int_\Lambda \sigma(x) v(x) \dx$. We use the cylindrical coordinates defined by $X$. Integrating along the radial direction, we have that: 
\[v(0,0,s) = v(r,\theta,s) - \int_0^r \frac{\partial v}{\partial r}(t,\theta,s)\dt.\] 
Hence: 
\begin{align*}
(\sigma\delta_\Lambda)(v) &=  \int_0^L \sigma(s) \left[ v(r,\theta,s) - \int_0^r \frac{\partial v}{\partial r}(t,\theta,s)\dt\right] \ds \\ &\le \|\sigma\|_{L^2(\Lambda)}\left(\int_0^L \left[ v(r,\theta,s) - \int_0^r \frac{\partial v}{\partial r}(t,\theta,s)\dt\right]^2 \ds \right)^\frac{1}{2} \\ 
			&\le \sqrt{2} \|\sigma\|_{L^2(\Lambda)}\left(\int_0^L v(r,\theta,s)^2 \ds + \int_0^L \left[\int_0^r \frac{\partial v}{\partial r}(t,\theta,s)\dt\right]^2 \ds \right)^\frac{1}{2} \end{align*} 
Now, we square this expression and we integrate in $C(\Lambda,R)$ for some $R\le R_0$, obtaining: 
\begin{align*}
	\pi R^2 (\sigma\delta_\Lambda)(v)^2 &\le 2\|\sigma\|_{L^2(\Lambda)}^2 \left(\int_0^{2\pi} \int_0^R\int_0^L v(r,\theta,s)^2 \ds\, r\; \dr\dtheta \right.\\ 
	&\quad\quad\quad\quad\quad\quad + \left.\int_0^{2\pi} \int_0^R \int_0^L \left[\int_0^r \frac{\partial v}{\partial r}(t,\theta,s) \dt \right]^2\ds \, r\; \dr \dtheta\right) \\ 
	&\le 2\|\sigma\|_{L^2(\Lambda)}^2\left(\|v\|_{L^2(C(\Lambda,R))}^2 +\int_0^{2\pi} \int_0^R \int_0^L\left[\int_0^r \frac{\partial v}{\partial r}(t,\theta,s) \dt\right]^2\ds\, r\; \dr \dtheta\right) 
\end{align*} 
We apply inequality \eqref{Hardy} with $p=q=2$, $\omega(t)=t$, $\nu(t) = t^{1-2\gamma}$ to the second term on the right-hand side, obtaining: 
\[\int_0^R \left[\int_0^r \frac{\partial v}{\partial r}(t,\theta,s)
\dt\right]^2 r \dr \le C(\gamma,R)^2 \int_0^R \left|\frac{\partial v}{\partial r}(t,\theta,s)\right|^2 r^{1-2\gamma}\dr,\] 
where $C(\gamma,R)\le k(2,2)D$. It is easy to check that $D=D(\gamma,R)=R^{1+\gamma}\frac{\gamma^{(\gamma-1)/2}}{2(\gamma+1)^{(\gamma+1)/2}}$, and $k(2,2)=2$.  Thus: 
\[C(\gamma,R) \le R^{1+\gamma}\frac{\gamma^{(\gamma-1)/2}}{(\gamma+1)^{(\gamma+1)/2}}.\]
Applying this inequality in the estimate above we have: 
\[|(\sigma\delta_\Lambda)(v)| \le \sqrt{\frac{2}{\pi}}\|\sigma\|_{L^2(\Lambda)}
\left(R^{-2}\|v\|^2_{L^2(C(\Lambda,R))} + R^{2\gamma}c(\gamma)\|\nabla
v\|^2_{L^2_{-\gamma}(C(\Lambda,R))} \right)^\frac{1}{2},\] 
with $c(\gamma)=\gamma^{\gamma-1}/(\gamma+1)^{\gamma+1}$, and
the result is proven. 

It is important to notice that $c(\gamma)\to\infty$ as
$\gamma\to 0^+$, so some weight is needed for the estimate
to hold. \end{proof}
 
The well-posedness of the weak problem is a direct consequence of the previous theorem:

\begin{thm} Let $\Omega\subset \R^3$ a $\mathcal{C}^{1,1}$ domain or a convex polyhedron and $0<\gamma<1$. Then problem \eqref{weakproblem} admits a unique solution $u\in	H^1_{\gamma}(\Omega)$ satisfying:
	\[\|u\|_{H^1_{\gamma}(\Omega)}\le C\|\sigma\|_{L^2(\Lambda)},\]
	where $C$ depends on $\gamma$ and $R_0$, and tends to $\infty$ as $\gamma\to 0$. 
\end{thm} 
	\begin{proof} The result follows directly from \cite[Corollary 2.7]{DDO} (see also \cite[Theorem 2.8]{DDO}) and  Theorem \ref{theorem continuity delta} above. Indeed, in the particular case $p=2$ and $\omega=d^\gamma$,  \cite[Corollary 2.7]{DDO} establishes that the problem $-\Delta u = \mu$ with $\mu\in (H^1_{-\gamma}(\Omega))'$ admits a unique solution in $H^1_\gamma(\Omega)$ satisfying the a priori estimate:
\[\|u\|_{H^1_\gamma(\Omega)}\le C\|\mu\|_{(H^1_{-\gamma}(\Omega))'}.\]
On the other hand, Theorem \ref{theorem continuity delta} shows that
\[\|\sigma\delta_\Lambda\|_{(H^1_{-\gamma}(\Omega))'}\le C \|\sigma\|_{L^2(\Omega)},\] 
with $C$ depending on the radius $R_0$. Taking $\mu=\sigma\delta_\Lambda$, and combining these results we obtain the theorem. 
\end{proof}

\section{Approximating problem}\label{section approximating problem}  

Our approach is based on the study of a regularized version of problem \eqref{problem}. We consider
the function $\phi_n\in C^\infty_0(\R^n)$: \[\phi_n(x) = \left\{
\begin{array}{cl} C e^{\frac{1}{|x|-1}} & |x|<1 \\ 0
& |x|\ge 1 \end{array},\right.\] where the constant $C=C(n)$
is chosen so that $\int\phi_n = 1$. Then, $\phi_{n,\rho}(x) =
\rho^{-n}\phi_n(x/\rho)$ is an approximation of the $n$-dimensional Dirac
delta supported at the origin that works under convolution as
a mollifier of well known properties (see for example
\cite[Appendix C.4]{Evans_PDE}). An important and easy to check property of $\phi_{n,\rho}$ is that: 
\begin{equation}\label{estimate phi}
	|D^\beta \phi_{n,\rho}|\le C\rho^{-n-|\beta|}.
\end{equation}

We take $\phi_{1,\rho}$ and
$\phi_{2,\rho}$ two versions onf $\phi_{n,\rho}$ for $n=1$ and
$n=2$ respectively. Then,  consider an approximation
of $\sigma\delta_\Lambda$ that we define in terms of the cylindrical
coordinates given by $X$: 
\[\sigma_\rho(r,\theta,s) =
\phi_{2,\rho}(r)\int_\rho^{L-\rho}\sigma(t)\phi_{1,\rho}(s-t)\dt.\] 
It is clear that supp$(\sigma_\rho)=C(\Lambda,\rho)$. Moreover, the integral factor is a convolution along the $s$ axis, whereas for each $s$,
$\phi_{2,\rho}$ is an approximation of a two-dimensional Dirac delta on
the plane of versors ${\bf n}(s)$ and ${\bf b}(s)$.  In order to
be able to use cylindrical coordinates, we assume that
$\rho<R_0$. 

In the sequel we will use extensively that the domain of integration of the integral in $\sigma_\rho$ is narrowed by the support of $\phi_{1,\rho}$. Indeed, for every fixed $s$, supp$(\phi_{1,\rho}(s-\cdot)) = (s-\rho,s+\rho)$, so we define: 
\[I_\rho(s) = (s-\rho,s+\rho).\]
When taking norms of $\sigma_\rho$, we will apply many times  Fubini's Lemma to two integrals along the $s$ axis. For this, it is useful to observe that: 
\[\{s:\,t\in I_\rho(s)\} = I_\rho(t).\]
An important fact is that for every $s$ and every $t$, $|I_\rho(s)|\sim \rho$. 

When necessary, we assume that $\sigma$ is extended by zero outside of the interval $[0,L]$. 

The following lemma proves that $\sigma_\rho$ is indeed an approximation of $\sigma\delta_\Lambda$.

\begin{lemma}\label{lemma sigma_rho approx} Let $v\in H^1_{-\gamma}(\Omega)$,
then: \[\lim_{\rho\to 0^+}\int \sigma_\rho(x) v(x) \dx =
(\sigma\delta_\Lambda)(v).\] 
\end{lemma} 
\begin{proof} By a density
argument, it is enough to consider $v\in C^\infty(\Omega)$.
When integrating only along the $s$ axis, we simplify the
notation writing $v(s)$ instead of $v(0,0,s)$. Integrating in
cylindrical coordinates, we have: \begin{align*} \int
	\sigma_\rho(x)v(x)\dx &= \int_0^{2\pi}\int_0^\rho \int_0^L
	r\phi_{2,\rho}(r)\int_\rho^{L-\rho}\sigma(t)\phi_{1,\rho}(s-t)\,\dt	v(r,\theta,s) \,\ds\dr\dtheta  \\ 
		&= \int_0^{2\pi}\int_0^\rho	\int_0^L r\phi_{2,\rho}(r)\int_\rho^{L-\rho}\sigma(t)\phi_{1,\rho}(s-t)\,\dt \Big[v(r,\theta,s)-v(s)\Big]\, \ds\dr\dtheta \\
	&\quad\quad\quad + \int_0^{2\pi}\int_0^\rho \int_0^L r\phi_{2,\rho}(r)\int_\rho^{L-\rho}\sigma(t)\phi_{1,\rho}(s-t)\,\dt \, v(s) \, \ds\dr\dtheta \\ 
	&= \int_0^{2\pi}\int_0^\rho\int_0^L r\phi_{2,\rho}(r) \int_\rho^{L-\rho}\sigma(t)\phi_{1,\rho}(s-t)\,\dt \Big[v(r,\theta,s)-v(s)\Big] \,\ds\dr\dtheta \\
	&\quad\quad\quad + \int_0^L \int_\rho^{L-\rho}\sigma(t)\phi_{1,\rho}(s-t)\dt 	v(s) 	\,\ds =: 	I + 	II,
\end{align*} where in the
last step we used that $\phi_{2,\rho}$ integrates $1$. We begin by proving that $I\to 0$ as $\rho\to 0$. As in Theorem \ref{theorem continuity delta}, we use that: \[v(r,\theta,s)-v(s) = \int_0^r \frac{\partial v}{\partial r}(s,\xi,\theta) \dxi.\]

Taking into account that $|\phi_{2,\rho}|\le C\rho^{-2}$ and applying the Cauchy-Schwartz
inequality, we have that 
\begin{align*}
	I &\le C\rho^{-2} \int_0^{2\pi} \int_0^L\int_\rho^{L-\rho}\sigma(t)\phi_{1,\rho}(s-t)\,\dt \ds \int_0^\rho \int_0^r 	\frac{\partial v}{\partial r}(s,\xi,\theta) \,\dxi\, r\;\dr \dtheta \\  
		&\le C\rho^{-2} \int_0^{2\pi}	\int_0^L 	\int_\rho^{L-\rho}\sigma(t)\phi_{1,\rho}(s-t)\,\dt\ds \left(\int_0^\rho 	\bigg|\int_0^r \frac{\partial v}{\partial r}(s,\xi,\theta) \,\dxi\bigg|^2	r \,\dr\right)^\frac{1}{2}\left(\int_0^\rho 	r\,\dr\right)^\frac{1}{2}\dtheta \\ 
		&\le C\rho^{-1} \int_0^{2\pi} \int_0^L \int_\rho^{L-\rho}\sigma(t)\phi_{1,\rho}(s-t)\,\dt\ds \left(\int_0^\rho 	\bigg|\int_0^r\frac{\partial v}{\partial r}(s,\xi,\theta)\,\dxi\bigg|^2 r\,\dr\right)^\frac{1}{2} \dtheta.
\end{align*} 
We apply the Hardy inequality \eqref{Hardy} with $p=q=2$, $\omega(t)=t$,
$\nu(t)=t^{1-2\gamma}$, recalling that $C(\gamma)\le c(\gamma)\rho^{1+\gamma}$:
\begin{align*} 
	I &\le C\rho^{\gamma} \int_0^{2\pi} \int_0^L \int_\rho^{L-\rho}\sigma(t)\phi_{1,\rho}(s-t)\,\dt\ds \left(\int_0^\rho \bigg| \frac{\partial v}{\partial r}(r,\theta,s)\bigg|^2 r^{1-2\gamma}\,	\dr\right)^\frac{1}{2}\,\dtheta. 
\end{align*}
	We continue by recalling the definition of $I_\rho(s)$ and applying
	once again the Cauchy-Schwartz inequality:
	\begin{align*} 
		I &\le C\rho^{\gamma} \left(\int_0^{2\pi}\int_0^L		\int_{I_\rho(s)\cap [\rho,L-\rho]} \Big|\sigma(t)\phi_{1,\rho}(s-t)\Big|^2\,\dt\ds\dtheta\right)^\frac{1}{2}	\\
		&\quad\quad\quad\quad\quad\cdot \left(\int_0^{2\pi}\int_0^L \int_{I_\rho(s)\cap[\rho,L-\rho]}\int_0^\rho \bigg|\frac{\partial v}{\partial r}(r,\theta,s)\bigg|^2 r^{1-2\gamma} \,\dr\dt\ds\dtheta\right)^\frac{1}{2} \\ 
		&\le	C\rho^{\frac{1}{2}+\gamma}\|\nabla v\|_{L^2_{-\gamma}(C(\Lambda,\rho))}\left(\int_0^L\int_{I_\rho(s)\cap[\rho,L-\rho]}	\Big|\sigma(t)\phi_{1,\rho}(s-t)\Big|^2\,\dt\ds \right)^\frac{1}{2}
\end{align*} 
Finally, 	let us  apply Fubini's 	lemma and the estimate $|\phi_{1,\rho}|\sim \rho^{-1}$.
		\begin{align*} 
			I &\le C\rho^{\frac{1}{2}+\gamma} \|\nabla v\|_{L^2_{-\gamma}(\Omega)}\left(\int_0^L \sigma(t)^2 \int_{I_\rho(t)}
			\phi_{1,\rho}(s-t)^2\, \ds \dt \right)^\frac{1}{2} \\
			&\le C\rho^{\frac{1}{2}+\gamma} \|\nabla v\|_{L^2_{-\gamma}(\rho)}\left(\int_0^L \sigma(t)^2 \rho \rho^{-2}\,	\dt \right)^\frac{1}{2}\\ 
			&\le C\rho^{\gamma} \|\nabla v\|_{L^2_{-\gamma}(\rho)}\|\sigma\|_{L^2(\Lambda)},
		\end{align*} 
		and $I\to 0$ as $\rho\to 0$ for every $\gamma>0$. 
		
		On the other hand, it is easy to check that $II$ tends to the desired limit. Indeed, applying Fubini's lemma and the Cauchy-Schwartz inequality:

		\begin{align*}
		  \Big|II-\int_\Lambda \sigma v\Big| &= \Big|\int_0^L\sigma(t)v(t)\dt\Big| \\ 
			&= \Big|\int_0^L\int_{\rho}^{L-\rho}\sigma(t)\phi_{1,\rho}(s-t)\dt \,v(s)\,\ds - \int_0^L\sigma(t)v(t)\dt\Big| \\
			&\le \int_\rho^{L-\rho}|\sigma(t)|\bigg|\int_0^L\phi_{1,\rho}(s-t)v(s)\,\ds -v(t)\bigg|\,\dt + L\int_{[0,\rho]\cup[L-\rho,L]}|\sigma(t)v(t)|\,\dt \\
			&\le \|\sigma\|_{L^2(\Lambda)}\|\phi_{1,\rho}*v-v\|_{L^2(\Lambda)} + L\int_{[0,\rho]\cup[L-\rho,L]}|\sigma(t)v(t)|\,\dt.
		\end{align*} 
		The second term vanishes as its domain of integration does as $\rho\to 0$, whereas the first one vanishes thanks to well known properties of the
		convolution with mollifiers (see \cite[Appendix C.4]{Evans_PDE}). 
	\end{proof}
	
We consider the approximating problem:
\begin{equation}\label{approx problem} \left\{
			\begin{array}{cl}
				-\Delta u_\rho = \sigma_\rho &  \textrm{ in }\Omega,\\ u_\rho = 0 &
				\textrm{ in }\partial\Omega,
			\end{array}\right. 
\end{equation}
			
Since $\sigma_\rho\in C^\infty$,  problem \eqref{approx problem} has a unique solution $u_\rho\in C^\infty$. In the following section we study weighted norms of $u_\rho$ and its derivatives, with weights of the form $d^{2\gamma}$ and choose the exponent $\gamma$ so that we can take limit with $\rho\to 0$. 
		
\section{Isotropic regularity}\label{section isotropic}

Our main isotropic result is stated in terms of the Kondratiev-type spaces $K^m_\gamma(\Omega)$, defined as
\[K^m_\gamma(\Omega) := \{v:\Omega\to\R:\, d^{\gamma+|\beta|} D^\beta v \in L^2(\Omega), \; \forall \beta:\, |\beta|\le m\},\] 
equipped with the norm: 
\[\|v\|_{K^m_\gamma(\Omega)}^2 := \sum_{0\le |\beta|\le m} \int_\Omega |D^\beta v(x)|^2 d(x)^{2(\gamma+|\beta|)}\,\dx\]

We prove the following theorem:
\begin{thm}\label{theorem isotropic}
	If $\sigma \in L^2(\Lambda)$ and $\Omega$ is a domain of class $\mathcal{C}^{m-1,1}$, then $u_\rho\in K^m_\gamma(\Omega)$ for every $\gamma>-1$. Moreover, the following estimate holds
	\[\|u_\rho\|_{K^m_\gamma(\Omega)} \le C\|\sigma\|_{L^2(\Lambda)},\]
	with a constant $C$ independent of $\rho$. Therefore, the solution $u$ of the singular problem \eqref{problem} also belongs to $K^m_\gamma(\Omega)$. 

	Furthermore, the result is also true for $m=2$ if $\Omega$ is a convex polyhedron. 
\end{thm}

The rest of this section is devoted to the proof of this result, which is done through a series of lemmas. We begin by decomposing $u_\rho$ into two parts. 

It is well known, (see, for example \cite[Theorem 1.1]{GW_Green} and \cite[Section 2.4]{GT_PDEs}) that under very general assumptions on the domain $\Omega$, problem \eqref{approx problem} admits a Green function, $G:\Omega\times\Omega\to \R$ such that: 
\begin{equation}\label{u rep}
				u_\rho(x) = \int_\Omega G(x,y)\sigma_\rho(y)\dy. 
\end{equation}
				
Morever we have that 
	\[G(x,y) = \Gamma(x-y)+h(x,y),\] 
where, $\Gamma$ is the fundamental solution:
	\[\Gamma(x-y) = -\frac{1}{4\pi |x-y|},\] 
and $h(x,y)$ is a harmonic function satisfying the boundary condition  $h(x,y) = -\Gamma(x-y)$ for every fixed $y$:
	\[\left\{\begin{array}{rclc}
		\Delta_x h(x,y) &=& 0 & x\in \Omega\\ 
		h(x,y) &=& -\Gamma(x-y) &x\in \partial\Omega\end{array}\right.\]
Hence, we can separate the solution $u_\rho$ into two parts:
	\[u_\rho(x) = \int_\Omega \Gamma(x-y) \sigma_\rho(y)\dy + \int_\Omega h(x,y) \sigma_\rho(y)\dy =: u_{\rho}^\circ(x) + u_\rho^\partial (x). \]
					
The first part ($u_\rho^\circ$) satisfies $-\Delta u_\rho^\circ(x) = \sigma_\rho(x)$, whereas the second part ($u_\rho^\partial$) corrects the boundary values of $u_\rho^\circ$. In particular, taking into account the support of $\sigma_\rho$, we have that:
	\[u_\rho^\partial(x) = 	\int_{C(\Lambda,\rho)} h(x,y) \sigma_\rho(y) \dy.\]
The existence of $h(x,y)$, and consequently that of $G$, is guaranteed if every point in the boundary of $\Omega$ is a regular point. A classical result says that if $x\in \partial\Omega$ is the vertex of an open truncated cone contained in $\Omega^c$, then $x$ is regular \cite[Theorem 8.27]{Helms}. However additional regularity on the domain is necessary in order to control the norm of the derivatives of $u^\partial_\rho$. Lemma \ref{u_rho near boundary} provides such estimates. But first, let us prove an auxiliary lemmas that will be important thoughout the paper: 

\begin{lemma}\label{norm of sigma_rho}
	Let $\rho<R_0$, $\eta>-1$ and $\beta$ a multiindex with $|\beta|=k$. Then:
	\[\|D^{\beta}\sigma_\rho\|_{L^2_\eta(C(\Lambda,\rho))}\le C\rho^{-1-k+\eta} \|\sigma\|_{L^2(\Lambda)},\]
	and 
  \[\|\sigma_\rho\|_{L^1(C(\Lambda,\rho))}\le C \|\sigma\|_{L^2(\Lambda)}.\]
\end{lemma}
\begin{proof}
	\eqref{estimate phi} implies that $|D^\beta(\phi_{1,\rho}\phi_{2,\rho})|\le C\rho^{-3-k}$. Integrating in cylindrical coordinates, applying this estimate, the Cauchy-Schwartz inequality and Fubini's lemma we obtain: 
	\begin{align*}
		\|D^{\beta}\sigma_\rho\|_{L^2_\eta(C(\Lambda,\rho))}^2 
	 &= \int_0^\rho \int_0^{2\pi}\int_0^L \Big|\int_{\rho}^{L-\rho}\sigma(t)D^{\beta}(\phi_{1,\rho}(s-t)\phi_{2,\rho}(r))\,\dt\Big|^2 r^{1+2\eta}\,\ds\dtheta\dr \\
	&\le C\rho^{-6-2k} \int_0^\rho r^{1+2\eta}\,\dr \int_0^L\Big|\int_{I_\rho(s)\cap[\rho,L-\rho]}\sigma(t)\,\dt\Big|^2\,\ds \\
	&\le C\rho^{-6-2k} \rho^{2+2\eta} \int_0^L\int_{I_\rho(s)}|\sigma(t)|^2\,\dt|I_\rho(s)|\,\ds \\ 
	&\le C\rho^{-4-2k+2\eta} \rho \int_0^L |\sigma(t)|^2\int_{I_\rho(t)}\ds\,\dt \\ 
	&\le C\rho^{-2-2k+2\eta}\|\sigma\|_{L^2(\Lambda)}^2, 
	\end{align*}
	which competes the proof of the first estimate. The restriction $\eta>-1$ is necessary for the integrability of $r^{1+2\eta}$.
	The second estimate follows from the first one with $k=0$ and $\eta=0$. Applying the Cauchy-Schwartz inequality we have 
	\begin{align*}
		\|\sigma_\rho\|_{L^1_\eta(C(\Lambda,\rho))} &\le \|\sigma_\rho\|_{L^2(C(\Lambda,\rho))}|C(\Lambda,\rho)|^\frac{1}{2} \\
						&\le C\|\sigma\|_{L^2(\Lambda)}.
	\end{align*}
\end{proof}

Now, we can prove estimates for the derivatives of $u^\partial_\rho$. 

\begin{lemma}\label{u_rho near boundary} 
	Let $\Omega$ be a domain of class $\mathcal{C}^{m-1,1}$ and $\beta$ a multiindex with $|\beta|=k\le m$. Then, taking $\gamma>-1$, the following estimate holds: 
		\[\|D^\beta u_\rho^\partial\|_{L^2_\gamma(\omega)}\le C \|\sigma\|_{L^2(\Lambda)}, \] 
	where the constant $C$ depends on $R_0$, $m$, the distance from $\Lambda$ to $\partial\Omega$ and on $\gamma$, but is independent of $\rho$. 
	
	The result is also true for $m=2$ if $\partial\Omega$ is a convex polyhedron. 
\end{lemma} 
\begin{proof}
	Let us begin by writing: 
	\[\|D^\beta u^\partial_\rho\|_{L^2_\gamma(\Omega)} \le 
	\|D^\beta u^\partial_\rho\|_{L^2_\gamma(B(\Lambda,R_0))}+\|D^\beta u^\partial_\rho\|_{L^2_\gamma(\Omega\setminus B(\Lambda,R_0))}:= I + II.\]
	We recall that $h$ is a $C^\infty$ function on $\Omega\times\Omega$ (see, for example \cite[Chapter 29]{Treves}). Hence, we have that there is a constant $K_0$ depending on $\beta$ and $R_0$ such that 
	\[|D^\beta h(x,y)|\le K_0,\]
for every $(x,y)\in B(\Lambda,R_0)\times B(\Lambda,R_0).$
Then, applying the second estimate in Lemma \ref{norm of sigma_rho}, we have: 
\begin{align*}
	I^2 &= \int_{B(\Lambda,R_0)}\bigg|\int_{C(\Lambda,\rho)} D^\beta h(x,y) \sigma_\rho(y)\,\dy\bigg|^2 d(x)^{2\gamma}\,\dx \\ 
			&\le K_0^2 \int_{B(\Lambda,R_0)} \bigg|\int_{C(\Lambda,\rho)}\sigma_\rho(y)\,\dy\bigg|^2 d(x)^{2\gamma}\,\dx \\
			&\le CK_0^2 \|\sigma\|_{L^2(\Lambda)}^2\int_{B(\Lambda,R_0)}d(x)^{2\gamma}\,\dx.
\end{align*}
Furthermore, the condition $\gamma>-1$ implies that the weight is integrable in $B(\Lambda,R_0)$, which concludes the estimate for $I$. 

For $II$, let us observe that in $\Omega\setminus B(\Lambda,R_0)$ $d(x)\le C=C(\gamma,R_0)$. Hence, we can drop the weight: 
\[II \le C\|D^\beta u^\partial_\rho\|_{L^2(\Omega\setminus B(\Lambda,R_0))}. \]
Now, we invoke \cite[Theorem 2.5.1.1]{Grisvard} which provides a priori estimates for harmonic functions in terms of its boundary data. In particular, we have that if $\Omega$ is of class $\mathcal{C}^{m-1,1}$ then 
	\begin{equation}\label{a priori estimate h}
	\|h(\cdot,y)\|_{H^{m}(\Omega)} \le C\|\Gamma(\cdot-y)\|_{H^{m-\frac{1}{2}}(\partial\Omega)}.
\end{equation}
Since $\Gamma(\cdot-y)$ is a $C^\infty$ function over $\partial\Omega$ for every fixed $y$, we can take the maximum of the right hand side of \eqref{a priori estimate h} with $y$ in the closure of $C(\Lambda,R_0)$ obtaining a constant $K_1$ that depends on $m$ and on $\Omega$ such that 
\[\|h(\cdot,y)\|_{H^m(\Omega)}\le K_1\]
for every $y\in C(\Lambda,R_0)$. With this, we continue by applying the Cauchy-Schwartz inequality, the first estimate in Lemma \ref{norm of sigma_rho} and Fubini's Lemma: 
\begin{align*}
	II^2 &\le C \int_{\Omega\setminus B(\Lambda,R_0)}\bigg|\int_{C(\Lambda,\rho)} D^\beta h(x,y) \sigma_\rho(y)\,\dy\bigg|^2\,\dx \\
			&\le C \int_{\Omega\setminus B(\Lambda,R_0)}\int_{C(\Lambda,\rho)} |D^\beta h(x,y)|^2 \,\dy \int_{C(\Lambda,\rho)}|\sigma_\rho(y)|^2\,\dy\;\dx \\
			&\le C \rho^{-2}\|\sigma\|_{L^2(\Lambda)}^2 \int_{C(\Lambda,\rho)}\int_{\Omega\setminus B(\Lambda,R_0)} |D^\beta h(x,y)|^2 \,\dx\dy \\
		 &\le CK_1^2 \rho^{-2}\|\sigma\|_{L^2(\Lambda)}^2 |C(\Lambda,\rho)| \\
		&\le CK_1^2 \|\sigma\|_{L^2(\Lambda)}^2,
\end{align*}
which concludes the proof for domains of class $\mathcal{C}^{m-1,1}$. 

	For convex polyhedra, it suffices to show that \eqref{a priori estimate h} holds for $m=2$. The rest of the proof is the same. Since $\Gamma(\cdot-y)$ is smooth on $\partial\Omega$ for every $y\in C(\Lambda,R_0)$, \cite[Theorem 2]{BDM_Traces} says that we can find a function $\bar{h}_y\in H^2(\Omega)$ such that $\bar{h}_y|_{\partial\Omega} = -\Gamma(x-y)$ (see also \cite[Theorem 5]{BG_Traces} where a similar result is obtained for general three-dimensional Lipschitz domains). Moreover, we have the estimate:
	\[\|\bar{h}_y\|_{H^2(\Omega)}\le C\|\Gamma(\cdot-y)\|_{H^{\frac{3}{2}}(\partial\Omega)}.\] Now, applying the results of \cite[Section 4.3.1]{MR_Polyhedral} we can find $\mathcal{H}_y$ the solution of the problem
	\[\left\{\begin{array}{rcll}
			\Delta \mathcal{H}_y &=& -\Delta \bar{h}_y &\textrm{ in }\Omega \\
			\mathcal{H} &=& 0 & \textrm{ in } \partial\Omega,
	\end{array}\right.\]
	which in turn satisfies an estimate of the form: 
	\[\|\mathcal{H}\|_{H^2(\Omega)}\le C\|\Delta\bar{h}_y\|_{L^2(\Omega)}.\]
	It is clear that $h(x,y) = \mathcal{H}_y(x) + \bar{h}_y(x)$, and combining the a priori estimates for $H_y$ and $\bar{h}_y$ we obtain \eqref{a priori estimate h} with $m=2$, completing the proof. 
	\end{proof}

	Observe that the previous result shows that for $\gamma>-1$, $u_\rho^\partial\in K^m_\gamma(\Omega)$, provided that $\Omega$ is of class $\mathcal{C}^{m-1,1}$ or a convex polyhedron when $m=2$.

It remains to estimate the norm of $u_\rho^\circ$, that 
captures the singularity of the data. Naturally, this is much more complicated. We begin by discussing some preliminary results that will be crucial in the sequel. 

The following well known result due to Sawyer and Wheeden is proven in \cite[Theorem 1]{SW_fracint}:

\begin{thm}\label{theorem SW}
For $0<\alpha<n$, let $I_\alpha(f)$ denote the fractional integral in $\R^n$ applied to the function $f$: 
\[I_{\alpha}(f)(x) = \int_{\R^n} |x-y|^{\alpha-n}f(y) dx.\]
Let $1<p\le q<\infty$. If for some $\tau>1$,
\begin{equation}\label{Apqalpha}
 |Q|^{\frac{\alpha}{n}+\frac{1}{q}-\frac{1}{p}}\bigg(\frac{1}{|Q|}\int_Q w(x)^{q\tau}dx\bigg)^{\frac{1}{q\tau}}\bigg(\frac{1}{|Q|}\int_Q v(x)^{-p'\tau} dx\bigg)^{\frac{1}{p'\tau}}\le C_\tau, \tag{$A_{p,q,\tau}^{\alpha}$}
\end{equation}
for every cube $Q\subset\R^n$, then the weighted inequality:
\begin{equation}\label{SWinequality}
\bigg(\int_{\R^n} |I_{\alpha}f(x)|^q w(x)^q dx \bigg)^{\frac{1}{q}} \le
C \bigg(\int_{\R^n} f(x)^p v(x)^p dx \bigg)^{\frac{1}{p}},
\end{equation}
holds for every $f\geq 0$. 

Moreover, in the case $p<q$, condition \eqref{lemma Apqalpha} can be simplified by taking $\tau=1$ and with this modification it is not only sufficient but also necessary for \eqref{SWinequality} to hold.
\end{thm}

We want to apply this theorem with weights given by powers of the distance to $\Lambda$. The following lemma gives conditions on the exponents for  \eqref{Apqalpha} to hold for such weights. The proof is very similar to the one in \cite[Lemma 3.3]{DLG_Stokes}. However, in that paper the authors considered weights given  by powers of the distance to a set $E$ which is in turn contained in an Ahlfors regular set. Here, we state the result directly in terms of the Assouad dimension of $E$. Since it is not essential for the rest of the paper, we difer the definition of the Assouad dimension, as well as the proof of the lemma to the Appendix. For our purposes it suffices to observe that:
\begin{itemize}
	\item the Assouad dimension of a smooth curve is $1$, 
	\item the Assouad dimension of an isolated point is $0$.
\end{itemize}

\begin{lemma}\label{lemma Apqalpha}
Let $E\subset \R^n$, $\dim_A(E)$ its Assouad dimension, $1<p\le q<\infty$, and $0<\alpha<n$, satisfying the additional restriction: 
\begin{equation}\label{Apqalphacondalpha}
\frac{n}{p}-\frac{n}{q}<\alpha.
\end{equation} 
Let also $\eta\in\R$ and
\begin{equation}\label{eta*}
	\eta^* = \eta+\frac{n}{p}-\frac{n}{q}-\alpha.
\end{equation}
If the following conditions are satisfied
\begin{eqnarray}
	\eta<\frac{n-\dim_A(E)}{q},\label{Apqalphacondeta} \\
	\eta^*>-\frac{n-\dim_A(E)}{p'},\label{Apqalphacondeta*}
 \end{eqnarray}
 then inequality \eqref{SWinequality} holds with weights $w=d(\cdot,E)^{-\eta}$, $v=d(\cdot,E)^{-\eta^*}$. 
\end{lemma}

The global argument for estimating the norms of $u_\rho^\circ$ and its derivatives is as follows. We have the representation formula:
\[u_\rho^\circ(x) =  \int_{C(\Lambda,\rho)} \Gamma(x-y)\sigma_\rho(y)\dy.\]
Thanks to well known properties of the convolution, this also gives us a represention for the derivatives of $u_\rho^\circ$, of the form: 
\[D^\beta u_\rho^\circ(x) =  \int_{C(\Lambda,\rho)} \Gamma(x-y)D^{\beta}\sigma_\rho(y)\dy.\] 
Moreover, 
\[\Big|\int\Gamma(x-y) D^\beta\sigma_\rho(y)\dy \Big|\sim |I_2(D^\beta\sigma_\rho)(x)|.\]

We study the weighted norm of $D^\beta u_\rho^\circ$ in several steps, according to an appropriate partition of the domain. First, we consider a neighbourhood of the support of $\sigma_\rho$ given by $B(\Lambda,2\rho)=C(\Lambda,2\rho)\cup B^+(0,2\rho)\cup B^+(L,2\rho)$. There, we use the previous representations and deal with the fractional integral involved by means of Lemma \ref{lemma Apqalpha}. In a second step, integrating with some care, we control the norm on $B(\Lambda,R_0)\setminus B(\Lambda,2\rho)$. Finally, the norm on $\Omega\setminus B(\Lambda,R_0)$ is easily estimated since the domain of integration is far from the singularity.

The following remark simplifies Lemma \ref{lemma Apqalpha}, focusing on the particular case that we use in this section.

\begin{rmk}\label{remark eta line}
	We will apply Lemma \ref{lemma Apqalpha} with $E=\Lambda$, $\alpha=2$ and $1<p<q=2$. First, observe that condition \eqref{Apqalphacondalpha} is fulfilled for every such $p$. Moreover, thanks to \eqref{eta*}, \eqref{Apqalphacondeta} and \eqref{Apqalphacondeta*} can be combined, obtaining restrictions stated only in terms of $\eta$. In particular, we have: 
	\begin{equation}\label{condition eta line}
		\frac{3}{2}-\frac{1}{p}<\eta<1
	\end{equation}
	which gives a feasible restriction for $\eta$ for every $p<2$.
\end{rmk}

\begin{lemma}\label{u_rho near}
	Let $\beta$ be a multiindex with $|\beta|=k$ and 
	\begin{equation}\label{gamma cond}
		\gamma>k-1.
	\end{equation}
		Then
	\begin{equation}\label{estimates for urhocirc}
		\|D^\beta u_\rho^\circ\|_{L^2_{\gamma}(C(\Lambda,2\rho))}\le C \|\sigma\|_{L^2(\Lambda)},
	\end{equation}
	where the constant $C$ is independent of $\rho$. 
\end{lemma}
\begin{proof}
	As mentioned above, we have that
	\[|D^\beta u_\rho^\circ(x)| \sim I_2(D^{\beta}\sigma_\rho)(x).\]

	Hence, by the dual characterization of the norm we have that
	\[\|D^\beta u_\rho^\circ\|_{L^2_\gamma(\rho)}\le 
	C\|I_2(D^{\beta}\sigma_\rho)\|_{L^2_\gamma(\rho)} = \sup_{g:\|g\|_{L^2_{-\gamma}(C(\Lambda,2\rho))}=1}\int_{C(\Lambda,2\rho)} g(x)I_2(D^{\beta}\sigma_\rho)(x)\dx.\]
	We continue by choosing 
	\[\eta=\frac{3}{2}-\frac{1}{p}+\varepsilon,\] 
	for some $\varepsilon>0$ small enough so that $\eta$ satisfies \eqref{condition eta line} for some $p<2$ to be determined later. Taking into account \eqref{eta*}, this gives:
	\[\eta^* = \varepsilon-\frac{2}{p'}.\]

 Then, applying Fubini's Lemma, multiplying by $d(x)^{\eta}d(x)^{-\eta}$ and using the H\"older inequality
	\begin{align*}
		\int_{C(\Lambda,2\rho)} g(x)I_2 (D^{\beta}\sigma_\rho)(x)\dx &=
		\int_{C(\Lambda,2\rho)} g(x)\int_{C(\Lambda,\rho)}|x-y|^{-1}D^{\beta}\sigma_\rho(y)\dy\dx \\
		&= \int_{C(\Lambda,\rho)} D^{\beta}\sigma_\rho(y)\int_{C(\Lambda,2\rho)}|x-y|^{-1}g(x)\dx\dy \\
		&= \int_{C(\Lambda,\rho)} D^{\beta}\sigma_\rho(y) I_2(g)(y)\dy \\
		&\le \|D^{\beta}\sigma_\rho\|_{L^2_\eta(C(\Lambda,\rho))}  \|I_2(g\chi_{C(\Lambda,2\rho)})\|_{L^2_{-\eta}(C(\Lambda,\rho))} = \bigstar.
	\end{align*}
	The first factor can bounded by Lemma \ref{norm of sigma_rho} giving
	\begin{equation}\label{first factor}
		\|\sigma_\rho\|_{L_{\eta}^2(\Lambda)}\le C\rho^{\eta-1-k}\|\sigma\|_{L^2(\Lambda)} = C\rho^{\frac{1}{2}-\frac{1}{p}+\varepsilon-k}\|\sigma\|_{L^2(\Lambda)}
	\end{equation}
	On the other hand, for the norm of $I_2(g\chi_{C(\Lambda,2\rho)})$, we apply Lemma \ref{lemma Apqalpha} with $p<q=2$ and $\eta^*$ as above. Then, we apply the H\"older inequality with exponents $2/p$ and $2/(2-p)$. Thus, we obtain
	\begin{align*}
		\|I_2(g\chi_{C(\Lambda,2\rho)})\|_{L^2_{-\eta}(C(\Lambda,\rho))}&\le C\|g\|_{L^p_{-\eta^*}(C(\Lambda,2\rho))} \\
			&=\Big(\int_{C(\Lambda,2\rho)} |g(x)|^p d(x)^{-p\eta^*}\dx\Big)^\frac{1}{p} \\
			&= \Big(\int_{C(\Lambda,2\rho)} g(x)^p d(x)^{-p\gamma}d(x)^{p(\gamma-\eta^*)}\dx\Big)^\frac{1}{p} \\
			&\le \Big(\int_{C(\Lambda,2\rho)} g(x)^2 d(x)^{-2\gamma}\dx\Big)^\frac{1}{2}\Big(\int_{C(\Lambda,2\rho)} d(x)^{\frac{2p}{2-p}(\gamma-\eta^*)}\dx\Big)^{\frac{2-p}{2p}}\\
			&=\Big(\int_{C(\Lambda,2\rho)} d(x)^{\frac{2p}{2-p}(\gamma-\eta^*)}\dx\Big)^{\frac{2-p}{2p}}, 
	\end{align*}
	where in the last step we used that $\|g\|_{L^2_{-\gamma}(C(\Lambda,2\rho))}=1$. Now we integrate in cylindrical coordinates and replace $\eta^*$ by its value, obtaining 
	\begin{align*}
		\|I_2(g\chi_{C(\Lambda,2\rho)})\|_{L^2_{-\eta}(C(\Lambda,\rho))}&\le \bigg(\int_0^{2\pi}\int_0^L\int_0^{2\rho} r^{\frac{2p}{2-p}(\gamma-\varepsilon+\frac{2}{p'})+1}\dr\ds\dtheta\bigg)^{\frac{2-p}{2p}}.
	\end{align*}
	For the integral to be finite, we need the integrability condition: 
	\[\frac{2p}{2-p}\Big(\gamma-\varepsilon+\frac{2}{p'}\Big)+1>-1,\]
	which through some simple calculations is shown to be equivalent to: 
	\[\gamma>-1+\varepsilon.\]
	With this, we complete the integral, obtaining: 
	\[\|I_2(g\chi_{C(\Lambda,2\rho)})\|_{L^2_{-\eta}(C(\Lambda,\rho))}\le C\rho^{\gamma+1-\varepsilon}.\]
	
	Joining the estimates for both factors in $\bigstar$ we conclude: 
	\begin{align*}\int_{C(\Lambda,2\rho)} g(x)I_2 (D^{\beta}\sigma_\rho)(x)\dx &\le C\rho^{\frac{1}{2}-\frac{1}{p}+\varepsilon-k}\rho^{\gamma+1-\varepsilon}\|\sigma\|_{L^2(\Lambda)} \\
	&\le C\rho^{\gamma+1-k+\frac{1}{2}-\frac{1}{p}}\|\sigma\|_{L^2(\Lambda)}.
\end{align*}
In order to complete the estimate we need to prevent the constant from going to infinity as $\rho$ vanishes. For this, we take into account that $\gamma+1-k>0$, so we can choose $p<2$ close enough to $2$ such that the exponent of $\rho$ is nonnegative. This concludes the proof. Since the restriction $\gamma>-1+\varepsilon$ is needed for $\varepsilon>0$ arbitrarily small, the result holds for every $\gamma>-1$. 
\end{proof}

Using the same ideas we can estimate the norm of $u^\circ_\rho$ and its derivatives near the endpoints of $\Lambda$:

\begin{lemma}\label{u_rho near ball}
	Let $\beta$ be a multiindex with $|\beta|=k$ and $\gamma$ satisfying \eqref{gamma cond}. 	Then the following estimates hold:
	\begin{equation}
		\|D^\beta u_\rho^\circ\|_{L^2_\gamma(B^+_{2\rho}(0)\cup B^+_{2\rho}(L))}\le C \|\sigma\|_{L^2(\Lambda)}
	\end{equation}
\end{lemma}
\begin{proof}
	We only estimate the norm over $B^+(0,2\rho)$, the other part is completely analogous.

	The proof is the same as the one in the previous lemma. We apply the dual characterization of the norm and Fubini's Lemma, arriving at an anologous to $\bigstar$. The first factor is once again estimated by Lemma \ref{norm of sigma_rho} and the second one by Lemma \ref{lemma Apqalpha}. This leads to the estimate: 
	\[\|I_2(g\chi_{B^+(0,\rho)})\|_{L^2_{-\eta}(C(\Lambda,\rho))}\le \Big(\int_{B^+(0,2\rho)} d(x)^{\frac{2p}{2-p}(\gamma-\eta^*)}\dx\Big)^{\frac{2-p}{2p}}.\]
	The only variation with respect to Lemma \ref{u_rho near} is that the integral is taken over $B^+(0,2\rho)$ so spherical coordinates are needed instead of cylindrical ones. This does not modify the final result.  
	Indeed, integrating in spherical coordinates we obtain
	\[\|I_2(g\chi_{B^+(0,\rho)})\|_{L^2_{-\eta}(C(\Lambda,\rho))}\le C\rho^{3\frac{2-p}{2p}+\gamma-\eta^*}=C\rho^{\frac{3}{p}-\frac{3}{2}+\gamma-\varepsilon+\frac{2}{p'}}.\] 
	Joining this with the estimate for the first factor in $\bigstar$ and applying \eqref{eta*} we have:
	\begin{align*}
		\int_{C(\Lambda,2\rho)} g(x)I_\alpha (D^{\beta''}\sigma_\rho)(x)\dx 
		&\le C\rho^{\frac{1}{2}-\frac{1}{p}+\varepsilon-k} \rho^{\frac{3}{p}-\frac{3}{2}+\gamma-\varepsilon+\frac{2}{p'}}\|\sigma\|_{L^2(\Lambda)} \\
		&\le C\rho^{\gamma-1+k}\|\sigma\|_{L^2(\Lambda)}\\
		&\le C\|\sigma\|_{L^2(\Lambda)},
	\end{align*}
where in the last step we applied condition \eqref{gamma cond}. 
\end{proof}

\begin{lemma}\label{u_rho not so near}
	Let $\beta$ be a multiindex with $|\beta|=k$ and $\gamma$ satisfying \eqref{gamma cond}. 	Then the following estimate holds:
	\begin{equation}
		\|D^\beta u_\rho^\circ\|_{L^2_\gamma(B(\Lambda,R_0)\setminus B(\Lambda,2\rho))}\le C \|\sigma\|_{L^2(\rho)},
	\end{equation}
	with a constant $C$ independent of $\rho$. 
\end{lemma}
\begin{proof}
	Here, it is convenient to apply the derivatives to the kernel $\Gamma$, which yelds 
	\[D^\beta u_\rho^\circ(x) = \int_{C(\Lambda,\rho)}D^\beta\Gamma(x-y)\sigma_\rho(y)\dy \le C \int_{C(\Lambda,\rho)}\frac{1}{|x-y|^{1+k}}|\sigma_\rho(y)|\dy.\]
We consider only the case where $x\in C(\Lambda,R_0)\setminus C(\Lambda,2\rho)$. The estimates for $x\in B^+_{R_0}(0)\setminus B^+_{2\rho}(0)$ of $x\in B^+_{R_0}(L)\setminus B^+_{2\rho}(L)$ are obtained following the same arguments, but integrating in spherical coordinates instead of cylindrical ones. 

Since we have the variables $y\in C(\Lambda,\rho)$ and $x\in C(\Lambda,R_0)\setminus C(\Lambda,2\rho)$, let us denote $s_x$, $r_x$ and $\theta_x$ the tangential, radial and angular coordinates corresponding to $x$ and $s_y$, $r_y$ and $\theta_y$ the ones corresponding to $y$.

We begin by observing that:
\begin{equation}\label{est mod qrho}
	|\sigma_\rho(y)|\le  C\rho^{-3} \int_{I_\rho(s_y)} |\sigma(t)|\,\dt.
\end{equation}  

Now, we separate $C(\Lambda,\rho)$ into two parts, depending on $x$. Namely: 
\[C(\Lambda,\rho) = \{y\in C(\Lambda,\rho):\,|s_y-s_x|<r_x\}\cup\{y\in C(\Lambda,\rho):\,|s_y-s_x|\ge r_x\} =: A \cup B,\]
which leads us to
\[\int_{C(\Lambda,\rho)}\frac{1}{|x-y|^{1+k}}|\sigma_\rho(y)|\dy = \int_{A} \frac{1}{|x-y|^{1+k}}|\sigma_\rho(y)|\dy+ \int_{B} \frac{1}{|x-y|^{1+k}}|\sigma_\rho(y)|\dy =: I + II.\]

We need to estimate the weighted norm of $I$ and $II$. For $I$, observe that if $y\in A$, then
\begin{align*}
	|x-y| &\le |x-\lambda(s_x)|+|\lambda(s_x)-\lambda(s_y)|+|\lambda(s_y)-y|  
	\le  r_x + r_x + r_y \le 2r_x + \rho \le \frac{5}{2}r_x. 
\end{align*}
Moreover $|x-y|\ge r_x-\rho > r_x/2$, so we have that $|x-y|\sim r_x$ for every $y\in A$. Consequently, applying \eqref{est mod qrho} and integrating in cylindrical coordinates, we obtain
\begin{align*}
	I &\le C\rho^{-3}\int_A \frac{1}{|x-y|^{1+k}}\int_{I_\rho(s_y)}|\sigma(t)|\,\dt\dy \\
	  &\le C\rho^{-3}r_x^{-1-k}\int_A\int_{I_\rho(s_y)}|\sigma(t)|\,\dt\dy \\
		&\le C\rho^{-3}r_x^{-1-k}\int_{s_y:|s_y-s_x|<r_x}\int_0^\rho\int_0^{2\pi} \int_{I_\rho(s_y)}\sigma(t)\,\dt\;r_y \dtheta_y\dr_y\ds_y \\
		&\le C\rho^{-3}r_x^{-1-k}\rho^2\int_{s_y:|s_y-s_x|<r_x} \int_{I_\rho(s_y)}|\sigma(t)|\,\dt\ds_y \\
		&\le C\rho^{-1}r_x^{-1-k}\int_{t:|t-s_x|<r_x+\rho}|\sigma(t)| \int_{I_\rho(t)}\,\ds_y \dt \\
		&\le Cr_x^{-1-k}\int_{t:|t-s_x|<r_x+\rho}|\sigma(t)|\, \dt,
\end{align*} 
where in the last steps we applied Fubini's Lemma and used that $|I_\rho(t)| \sim \rho$ for every $t$. We conclude by applying the Cauchy-Schwartz inequality and the fact that $\rho<r_x$, which implies that $|\{t:|t-s_x|<r_x+\rho\}|\le C r_x$.
\begin{align*}
	I&\le Cr_x^{-1-k}\Big(\int_{t:|t-s_x|<r_x+\rho}\sigma(t)^2\, \dt\Big)^\frac{1}{2}|\{t:|t-s_x|<r_x+\rho\}|^\frac{1}{2} \\ 
	&\le Cr_x^{-\frac{1}{2}-k}\Big(\int_{t:|t-s_x|<r_x+\rho}\sigma(t)^2\, \dt\Big)^\frac{1}{2} \\ 
\end{align*}

Inserting this estimate in the norm, integrating in cylindrical coordinates and  applying Fubini's Lemma on the integrals along the $s$ axis, we have
\begin{align*}
	\|I\|_{L^2_\gamma(C(\Lambda,R_0)\setminus C(\Lambda,2\rho))}^2 &\le C\int_{C(\Lambda,R_0)\setminus C(\Lambda,2\rho)} r_x^{-1-2k}\int_{t:|t-s_x|<r_x+\rho}\sigma(t)^2\, \dt\, d(x)^{2\gamma} \,\dx \\ 
	&\le C\int_{0}^L\int_{2\rho}^{R_0}\int_0^{2\pi} \int_{t:|t-s_x|<r_x+\rho}\sigma(t)^2\, \dt\, r_x^{2\gamma-1-2k+1} \,\dtheta_x\dr_x\ds_x \\ 
  &\le C\int_0^{R_0} \int_{0}^L \sigma(t)^2 \int_{s_x:|t-s_x|<r_x+2\rho}\,\ds_x \dt\, r_x^{2\gamma-2k} \,\dr_x \\
	&\le C\int_0^{R_0} \int_{0}^L \sigma(t)^2 \dt\, r_x^{2\gamma-2k+1} \,\dr_x \\ 
	&\le C\|\sigma\|_{L^2(\Lambda)}^2,
\end{align*}
where in the last step we used the integrabilty condition \eqref{gamma cond}. 

For $II$ the idea is quite similar, but we need to further decompose $B$ into several sets: 
\[B = \bigcup_{j=0}^{J-1} \Big\{y\in C(\Lambda,\rho):\, 2^jr_x<|s_y-s_x|\le 2^{j+1}r_x\Big\} =: \bigcup_{j=0}^{J-1} B_j,\]
where $J$ is the minimum integer such that $2^J r_x \ge L$. 
If $x$ is near an extreme point of $\Lambda$, each $B_j$ is a cylinder around $\Lambda$ of height $\sim 2^j r_x$ and radius $\rho$, at a distance $\sim 2^j r_x$ from $x$. On the other hand, if $x$ is near the center of $\Lambda$ each $B_j$ is formed by two of such cylinders (one on each side of $x$). We also have that for $y\in B_j$, $|x-y|\sim 2^j r_x$. .  The estimate on each $B_j$ is a copy of the estimate on $A$, but with $2^j r_x$ instead of $r_x$: 

\begin{align*}
	II &\le \sum_{j=0}^{J-1} \int_{B_j} \frac{1}{|x-y|^{1+k}} |\sigma_\rho(y)|,\dy \\
	   &\le \sum_{j=0}^{J-1} \rho^{-3} (2^j r_x)^{-1-k}\int_{B_j}\int_{I_\rho(s_y)}|\sigma(t)| \,\dt\dy\\ 
		&\le \sum_{j=0}^{J-1} \rho^{-3} (2^j r_x)^{-1-k}\int_{s_y:2^jr_x<|s_y-s_x|\le 2^{j+1}r_x}\int_0^\rho\int_0^{2\pi}\int_{I_\rho(s_y)}|\sigma(t)|\,\dt \dtheta_y\dr_y\ds_y \\ 
		 &\le \rho^{-1}\sum_{j=0}^{J-1} (2^j r_x)^{-1-k}\int_{t:2^jr_x-\rho<|t-s_x|\le 2^{j+1}r_x+\rho}|\sigma(t)|\int_{I_\rho(t)}\, \ds_y \dt \\
		&\le \sum_{j=0}^{J-1} (2^j r_x)^{-1-k}\int_{t:2^jr_x-\rho<|t-s_x|\le 2^{j+1}r_x+\rho}|\sigma(t)|\,\dt. \\ 
\end{align*}
We continue by applying the Cauchy-Schwartz twice: first to the integral and then to the summation:
\begin{align*} 
		\quad &\le \sum_{j=0}^{J-1} (2^j r_x)^{-\frac{1}{2}-k}\Big(\int_{t:2^jr_x-\rho<|t-s_x|\le 2^{j+1}r_x+\rho}\sigma(t)^2\,\dt\Big)^\frac{1}{2} \\ 
		 &\le J^\frac{1}{2} \bigg(\sum_{j=0}^{J-1} (2^j r_x)^{-1-2k}\int_{t:2^jr_x-\rho<|t-s_x|\le 2^{j+1}r_x+\rho}\sigma(t)^2\,\dt\bigg)^\frac{1}{2}.
\end{align*}
Now we proceed as in the estimation of the norm of $I$. The only remarkable difference is the appearence of the factor $J\sim|\log(r_x)|$. 

\begin{align*}
	\|II&\|_{L^2_\gamma(C(\Lambda,R_0)\setminus C(\Lambda,2\rho))}^2 \\
	&\le C\int_{C(\Lambda,R_0)\setminus C(\Lambda,2\rho)} J\sum_{j=0}^{J-1} (2^j r_x)^{-1-2k}\int_{t:2^jr_x-\rho<|t-s_x|\le 2^{j+1}r_x+\rho}\sigma(t)^2\,\dt\, d(x)^{2\gamma}\,\dx \\ 
	&\le C \int_{2\rho}^{R_0}\int_{0}^L\int_0^{2\pi} J\sum_{j=0}^{J-1} (2^j r_x)^{-1-2k}\int_{t:2^jr_x-\rho<|t-s_x|\le 2^{j+1}r_x+\rho}\sigma(t)^2\,\dt\,  r_x^{2\gamma+1} \,\dtheta_x\dr_x\ds_x \\
&\le  C    \int_{0}^{R_0}J\sum_{j=0}^{J-1}2^{-j(1+2k)}\int_{0}^L\int_{t:2^jr_x-\rho<|t-s_x|\le 2^{j+1}r_x+\rho}\sigma(t)^2\,\dt \ds_x\, r_x^{2\gamma-2k} \,\dr_x \\
&\le C  \int_{0}^{R_0}  J\sum_{j=0}^{J-1}2^{-j(1+2k)}\int_{0}^L \sigma(t)^2 \int_{s_x:2^jr_x-2\rho<|t-s_x|\le 2^{j+1}r_x+2\rho}\,\ds_x \dt\, r_x^{2\gamma-2k} \,\dr_x 
\end{align*}
\begin{align*}
  \quad	&\le C \|\sigma\|_{L^2(\Lambda)}^2  \int_{0}^{R_0} J\sum_{j=0}^{J-1}2^{-j(1+2k)} (2^jr_x)\, r_x^{2\gamma-2k} \,\dr_x \\
	&\le C  \|\sigma\|_{L^2(\Lambda)}^2   \int_{2\rho}^{R_0}   J\sum_{j=0}^{J-1}2^{-j(2k)}\, r_x^{2\gamma-2k+1} \,\dr_x.   
\end{align*} 
The summation over $j$ is finite for every $k\ge 1$. On the other hand, if $k=0$ then the summation is $J$, so we can continue the estimate assuming this worst possible case, and recalling that $J\sim |\log(r_x)|$: 

\begin{align*}
	\|II\|_{L^2_\gamma(C(\Lambda,R_0)\setminus C(\Lambda,2\rho))}^2 	&\le C  \|\sigma\|_{L^2(\Lambda)}^2   \int_{0}^{R_0}   |\log(r_x)|^2\, r_x^{2\gamma-2k+1} \,\dr_x   \\
	&\le C\|\sigma\|_{L^2(\Lambda)}^2,
\end{align*}
where in the last step we used that thanks to condition \eqref{gamma cond} the integral is finite for every $\rho$.  This concludes the proof. 
\end{proof}

Finally, we can estimate the derivatives of $u_\rho^\circ$ far from $\Lambda$: 

\begin{lemma}\label{u_rho far}
	Given $\beta$ a multiindex with $|\beta|=k$, the following estimate holds for every $\gamma$
	\[\|D^\beta u_\rho^\circ\|_{L^2_\gamma(\Omega\setminus B(\Lambda,R_0))}\le C\|\sigma\|_{L^2(\Lambda)},\]
	where the constant $C$ is independent of $\rho$. 
\end{lemma}
\begin{proof}
	We write, as in the previous lemma:
	\begin{align*} 
	 |D^\beta u_\rho^\circ(x)| &= \bigg|\int_{C(\Lambda,\rho)}D^\beta\Gamma(x-y)\sigma_\rho(y)\dy\bigg|\le \int_{C(\Lambda,\rho)}|x-y|^{-1-k}|\sigma_\rho(y)|\dy.
	\end{align*} 
	Since $x\in \Omega\setminus B(\Lambda,R_0)$ and $y\in C(\Lambda,\rho)$, we have that $|x-y|^{-1-k}\le CR_0^{-1-k}$. This and the second estimate in Lemma \ref{norm of sigma_rho} give
	\[ |D^\beta u_\rho^\circ(x)| \le C\|\sigma\|_{L^2(\Lambda)},\]
	where the constant $C$ depends on $R_0$  and on $k$, but not on $\rho$. Hence
	\begin{align*}
		\|D^\beta u_\rho^\circ\|_{L^2_\gamma(\Omega\setminus B(\Lambda,R_0))} &\le C\|\sigma\|_{L^2(\Lambda)}\|d^\gamma\|_{L^2(\Omega\setminus B(\Lambda,\rho))},
	\end{align*}
	and the result follows from the integrability of the weight.  
\end{proof}

Thus, we have concluded the proof of Theorem \ref{theorem isotropic}, which is obtained by the combination of Lemmas \ref{u_rho near boundary}, \ref{u_rho near}, \ref{u_rho near ball}, \ref{u_rho not so near} and \ref{u_rho far}.

\section{Anisotropic regularity}\label{section anisotropic}

For simplicity, in this section we assume that $\Lambda$ is a straight segment, namely,
\[\Lambda = \{\lambda(s): \; 0\le s\le L, \;\textrm{ with }\lambda(s)=(0,0,s)\}.\]
This assumption is not really necessary and it is only introduced in order to simplify some calculations. In the next section we discuss the case of a general curved fracture, where essentially the same results can be proven. More importantly, we assume that
 \[\sigma\in H^m(\Lambda)\quad \quad \textrm{ for some } m\ge 1.\]
We also make extensive use of the compact embedding $H^1(\Lambda)\subset L^\infty(\Lambda)$, which gives the following estimate: 
\begin{equation}\label{embedding}
	|\sigma^{(\ell)}(t)| \le C\|\sigma^{(\ell)}\|_{H^{1}(\Lambda)}  \quad \quad\forall \ell\le m-1. 
\end{equation}

Our anisotropic estimates follow from the simple observation that near $\Lambda$, the derivatives of the solution of problem \eqref{problem} with respect to $s$ are smoother than its derivatives in any other direction. However, as we shall see, for the derivatives with respect to $s$ a singularity arises near the extreme points of $\Lambda$. Consequently, we define the anisotropic Kondratiev type spaces as follows:

\begin{defi}
	Given a multi-index $\beta$, we can distinguish the derivatives along $s$, that we denote $\beta_s$ from the derivatives with respect to the other variables, that we denote $\beta_\perp$. With this notation we have $\beta=\beta_s + \beta_\perp$. We also denote $k_s=|\beta_s|$ and $k_\perp = |\beta_\perp|$.

	Furthermore, let  
	\[d_e(x)=\min\{|x-\lambda(0)|,|x-\lambda(L)|\}\]
	be the distance to the extreme points of $\Lambda$. Then, given a domain $E$ we denote $K_{\gamma,\mu}^m(E)$ the Kondratiev-type space formed by the functions $v$ such that the following norm is finite: 
	
	\[\|v\|_{K_{\gamma,\mu}^m(E)}^2 := \sum_{\beta:|\beta|\le m} \int_E |D^\beta v(x)|^2 d(x)^{2(\gamma+|\beta_{\perp}|)}d_e(x)^{2(\mu+|\beta_s|)}\,\dx.\]

	We also denote $L^2_{\gamma,\mu}(E)=K^0_{\gamma,\mu}(E)$. In particular $L^2_{0,\mu}$ is the space with weight $d_e^\mu$. 
\end{defi}

It is important to take into account that the solution of problem \eqref{problem} is smooth far from $\Lambda$. Indeed, Lemmas \ref{u_rho near boundary} and \ref{u_rho far} imply that $u\in H^m(\Omega\setminus B(\Lambda,R_0))$ provided that $\Omega$ is of class $\mathcal{C}^{m-1,1}$ (or a convex polyhedron in the case $m=2$). Consequently, we only analyse the anisotropic behaviour of $u$ in $B(\Lambda,R_0)$. Our main result is the following: 

\begin{thm}\label{main theorem anisotropic}
	Let $\sigma\in H^m(\Lambda)$ for some $m\ge 1$, and $\Omega$ a domain of class $\mathcal{C}^{m-1,1}$. If
	\begin{equation}\label{gamma cond anisotropic}
		\gamma>k_\perp-1,
	\end{equation}
	\begin{equation}\label{mu cond anisotropic}
		\mu>k_s-\frac{1}{2},
	\end{equation}
	then $u_\rho\in K_{\gamma,\mu}^m(B(\Lambda,R_0)$. Moreover: 
	\[\|u_\rho\|_{K_{\gamma,\mu}^m(B(\Lambda,R_0)}\le C\|\sigma\|_{H^m(\Lambda)},\]
	where the constant $C$ is independent of $\rho$. 

	Therefore, the solution $u$ of the singular problem \eqref{problem} also belongs to $K_{\gamma,\mu}^m(B(\Lambda,R_0))$.

	The result also holds for convex polyhedra with $m=2$. 
\end{thm}

The main difficulty of the proof lies in the necessity of handling the two weights. For this, it is convenient to consider an appropriate decomposition of a cylinder surrounding $\Lambda$. Let us begin by introducing for every $0\le a<b\le L$, the notation
\[\Lambda[a,b]=\{\lambda(s):\; a\le s\le b\},\]
which represents a curve contained in $\Lambda$. We also denote $C(\Lambda[a,b],\rho)$ the cylinder around $\Lambda[a,b]$. 

For the sake of simplicity and without loss of generality, we can assume that values of $\rho$ are chosen such that there is an integer $J$ satisfying $2^J\rho = L$. We define, for every $0<\delta\le R_0$ 
\begin{equation}\label{decomp}
	{\bf C}^\delta_0 := C(\Lambda[0,\rho],\delta),\quad\quad {\bf C}^\delta_j := C(\Lambda[2^{j-1}\rho,2^{j}\rho],\delta),\quad \textrm{ for }j=1,\dots,J.
\end{equation}
In this way, we have that
\[C(\Lambda,\delta) = \bigcup_{j=0}^J {\bf C}_j^\delta, \quad \textrm{ and }\quad C(\Lambda[0,\dfrac{L}{2}],\delta) = \bigcup_{j=0}^{J-1} {\bf C}_j^\delta.\]
We also define expanded versions of ${\bf C}_j^\delta$: 
\[\overline{\bf C}^\delta_j := \cup_{i=j-1}^{j+1}{\bf C}^\delta_i.\]

The advantage of this decomposition is that $d_e(x)$ can be regarded as essentially constant over ${\bf C}_j^\delta$ for every $j\ge 1$, and consequently this weight can be pulled out of the norm. For studying the norm in a neighbourhood of an extreme point (such as ${\bf C}_0^\delta$), we take into accont the following remark.

\begin{rmk}\label{two weight trick}
	Let us consider a neighbourhood of $\lambda(0)$. There, we can integrate in spherical coordinates $(r,\theta,\xi)$ where $\xi$ is the cenital angle. In this case, we have that $d_e(x) = r$ and $\frac{d(x)}{d_e(x)} =\sin(\xi)$, hence a product of powers of $d(x)$ and $d_e(x)$ can be written as follows: 
	\[d(x)^ad_e(x)^b = r^{a+b}\sin(\xi)^a.\]
\end{rmk}

Let us begin the proof of Theorem \ref{main theorem anisotropic} by obtaining an analogous to Lemma \ref{u_rho near boundary}

\begin{lemma}\label{u_rho near boundary anisotropic}
	Let $\Omega$ be a domain of class $\mathcal{C}^{m-1,1}$ and $\beta$ a multiindex with $|\beta|=k\le m$. Then, taking $\gamma>-1$ and $\mu>-\frac{1}{2}$, the following estimate holds: 
	\[\|D^\beta u_\rho^\partial\|_{L^2_{\gamma,\mu}(\Omega)}\le C \|\sigma\|_{L^2(\Lambda)}, \] 
	where the constant $C$ depends on $R_0$, $m$, the distance from $\Lambda$ to $\partial\Omega$ and on $\gamma$, but is independent of $\rho$. 
	
	The result is also true for $m=2$ if $\partial\Omega$ is a convex polyhedron. 
\end{lemma}
\begin{proof}
	The proof is completely analogous to the one of Lemma \ref{u_rho near boundary}. The only difference lies in the estimate of the term $I$, where it is necessary to prove the integrability of the weights in $B(\Lambda,R_0)$. For this let us split the integral into three subdomains: $B(\lambda(0),R_0)$, $B(\lambda(L),R_0)$ and $B(\Lambda,R_0)\setminus (B(\lambda(0),R_0)\cup B(\lambda(L),R_0))$. In the third one we have that $R_0\le d_e(x)\le {\rm diam}(\Omega)$, so $d_e^{2\mu}\le C$ and $d(x)^{2\gamma}$ is integrable for $\gamma>-1$. On the other hand, integrating in spherical coordinates and recalling Remark \ref{two weight trick} we have that: 
	\begin{align*}
		\int_{B(\lambda(0),R_0)}d(x)^{2\gamma}d_e(x)^{2\mu}\,\dx &= \int_0^{R_0}\int_0^{2\pi}\int_0^\pi r^{2\gamma+2\mu+2} \sin(\xi)^{2\mu+1}\,\dxi\dtheta\dr \\
&\le C \int_0^{R_0} r^{2\gamma+2\mu+2}\,\dr \int_0^\pi \sin(\xi)^{2\mu+1}\,\dxi, 
	\end{align*}
	and both integrals are finite under the conditions $\gamma>-1$ and $\gamma+\mu>-\frac{3}{2}$. The integral in $B(\lambda(L),R_0)$ can be estimated in the same way. 
\end{proof}

As in the isotropic case, the estimates for $u_\rho^\circ$ are much more complicated. We prove them in several lemmas. We begin by stating some auxiliary results that will be helpful in the sequel.

\begin{lemma}
	Let $\Lambda$ be the segment defined above, $\beta_s$ a multiindex corresponging to derivatives only in $s$ and $k_s = |\beta_s|$. We denote $\sigma^{(k_s)}$ the derivative of order $k_s$ of $\sigma$ and \[(\sigma^{(k_s)})_\rho = \int_\rho^{L-\rho}\sigma^{(k_s)}(t)\phi_{1,\rho}(s-t)\,\dt\;\phi_{2,\rho}(r)\]
the regularization of $\sigma^{(k_s)}$.
	We define, for $k_s\ge 1$ and $e\in\{\rho,L-\rho\}$ the following functions in cylindrical coordinates:
	\[\mathcal{E}^{k_s}_e(r,\theta,s)=\sum_{\ell=0}^{k_s-1} \sigma^{(\ell)}(e)\phi_{1,\rho}^{(k_s-1-\ell)}(s-e)\phi_{2,\rho}(r).\] 
	Then, if $\sigma\in H^m(\Lambda)$ and $k_s\le m$, the following identity holds for every $(r,\theta,s)\in C(\Lambda,R_0)$: 
	\begin{equation}\label{anisotropic deriv}
		D^{\beta_s} \sigma_\rho(r,\theta,s)  = (\sigma^{(k_s)})_\rho+ \mathcal{E}^{k_s}_\rho(r,\theta,s) - \mathcal{E}^{k_s}_{L-\rho}(r,\theta,s).
	\end{equation} 
\end{lemma}
\begin{proof}
	Let us begin considering the first order derivative with respect to $s$, i.e.: $k_s=1$. It is immediate that: 
	\[D^{\beta_s} \sigma_\rho(r,\theta,s) = \phi_{2,\rho}(r) \int_\rho^{L-\rho} \sigma(t)\phi_{1,\rho}'(s-t)\,\dt.\]
	Integrating by parts, we obtain the desired result:
	\begin{align*}
		D^{\beta_s} &u_\rho(r,\theta,s) = \phi_{2,\rho}(r) \bigg\{\int_\rho^{L-\rho} \sigma'(t)\phi_{1,\rho}(s-t)\,\dt - \sigma(t)\varphi_{1,\rho}(s-t)\Big|_{t=\rho}^{t=L-\rho}\bigg\} \\
		&=\phi_{2,\rho}(r) \bigg\{\int_\rho^{L-\rho} \sigma'(t)\phi_{1,\rho}(s-t)\,\dt + \sigma(\rho)\varphi_{1,\rho}(s-\rho) - \sigma(L-\rho)\varphi_{1,\rho}(s-L+\rho)\bigg\}.
	\end{align*}
	A simple induction argument gives the identities for derivatives of order $k_s>1$. 
\end{proof}

The next two lemmas are analogous to Lemma \ref{norm of sigma_rho} for the anisotropic case.

\begin{lemma}\label{norm sigma_rho anisotropic}
	Assuming $\sigma\in H^m(\Lambda)$ for some $m\ge 1$, let $0\le \ell<m$. We denote $\beta_\perp$ a multiindex corresponding only to derivatives in directions orthogonal to $\Lambda$, and $|\beta_\perp|=k_\perp$.  Then:
	\[|D^{\beta_\perp}(\sigma^{(\ell)})_\rho(r,\theta,s)| \le C\rho^{-2-k_\perp}\|\sigma\|_{H^{\ell+1}(\Lambda)}.\]
	As a consequence, we have that for every $\eta>-\frac{3}{2}$:
	\[\|D^{\beta_\perp}(\sigma^{(\ell)})_\rho\|_{L^2_{0,\eta}({\bf C}_0^\rho)} \le C\rho^{\eta-\frac{1}{2}-k_\perp}\|\sigma\|_{H^{\ell+1}(\Lambda)}.\] 
	Moreover, for $\eta>-\frac{1}{2}$:
	\[\|D^{\beta_\perp}(\sigma^{(\ell)})_\rho\|_{L^2_\eta({\bf C}_j^\rho)} \le C\rho^{\eta-1-k_\perp}(2^j\rho)^\frac{1}{2}\|\sigma\|_{H^{\ell+1}(\Lambda)}.\]
\end{lemma}
\begin{proof}
	The first estimate is derived from \eqref{estimate phi} and \eqref{embedding}: 
	\begin{align*}
		|D^{\beta_\perp}(\sigma^{(\ell)})_\rho(r,\theta,s)| &= \Big|D^{\beta_\perp}\phi_{2,\rho}(r)\int_{I_\rho(s)\cap[\rho,L-\rho]}\sigma^{(\ell)}(t)\phi_{1,\rho}(s-t)\,\dt\Big| \\
		&\le C\rho^{-3-k_\perp}\int_{I_\rho(s)}|\sigma^{(\ell)}(t)|\,\dt \\
		&\le C\rho^{-3-k_\perp} \|\sigma^{(\ell)}\|_{H^1(\Lambda)} |I_\rho(s)| \\
		&\le C\rho^{-2-k_\perp} \|\sigma\|_{H^{\ell+1}(\Lambda)}.
	\end{align*} 
  
	For the second estimate, we apply the first one. Then, we use that ${\bf C}_0^\rho\subset B(\lambda(0),\sqrt{2}\rho)$ and we integrate the weight in spherical coordinates:  
	\begin{align*}  
		\|D^{\beta_\perp}(\sigma^{(\ell)})_\rho\|_{L^2_{0,\eta}({\bf C}^\rho_0)} &\le C\rho^{-2-k_\perp}\|\sigma\|_{H^{\ell+1}(\Lambda)} \bigg(\int_{{\bf C}^\rho_0}d_e(x)^{2\eta} \,\dx \bigg)^\frac{1}{2}\\
			&\le C\rho^{-2-k_\perp} \|\sigma\|_{H^{\ell+1}(\Lambda)}\bigg(\int_0^{\sqrt{2}\rho} \int_0^\pi \int_0^{2\pi} r^{2\eta+2}\,\dtheta\dxi\dr\bigg)^\frac{1}{2} \\ 
	&\le C\rho^{-2-k_\perp}\|\sigma\|_{H^{\ell+1}(\Lambda)}\rho^{\eta+\frac{3}{2}} \\
	&\le C\rho^{\eta-\frac{1}{2}-k_\perp}\|\sigma\|_{H^{\ell+1}(\Lambda)}. 
	\end{align*}
	Naturally, the same estimate holds on the cylinder $C(\Lambda[L-\rho,L],\rho)$, which is a neighbourhood of the extreme point $\lambda(L)$. 

	Finally, for the third estimate, we apply once again the first one and we integrate in cylindrical coordinates: 
	\begin{align*}
		\|D^{\beta_\perp}(\sigma^{(\ell)})_\rho\|_{L^2_\eta({\bf C}_j^\rho)} &\le C\rho^{-2-k_\perp}\|\sigma\|_{H^{\ell+1}(\Lambda)}\bigg(\int_{{\bf C}_j^\rho}d(x)^{2\eta}\dx\bigg)^\frac{1}{2} \\
		&\le C\rho^{-2-k_\perp}\|\sigma\|_{H^{\ell+1}(\Lambda)}\bigg(2^j\rho\int_0^\rho r^{2\eta+1}\dr\bigg)^\frac{1}{2}\\
		&= C\rho^{\eta-1-k_\perp}(2^j\rho)^\frac{1}{2}\|\sigma\|_{H^{\ell+1}(\Lambda)}.
\end{align*}
\end{proof}

It is clear that in the second estimate ${\bf C}_0^\rho$ can be replaced by $\overline{\bf C}_0^\rho$ and in the third one ${\bf C}_j^\rho$ can be replaced by $\overline{\bf C}_j^\rho$ if $j\ge 1$.  

The terms $\mathcal{E}^{k_s}_e$ with $e=\rho$ and $e=L-\rho$ are symmetrical and can be treated in the same way. Hence, we establish our results only in terms of $\mathcal{E}^{k_s}_\rho$. 

\begin{lemma}\label{E estimate}
	Given $\beta_\perp$ a multiindex representing derivatives in directions orthogonal to $\Lambda$ with $|\beta_\perp|=k_\perp$ and $k_s\ge 1$, the following estimate holds: 
	\[|D^{\beta_\perp}\mathcal{E}^{k_s}_\rho(r,\theta,s)|\le C\rho^{-2-k_\perp-k_s} \|\sigma\|_{H^{k_s}(\Lambda)} \chi_{\overline{\bf C}_0^\rho}.\]
	As a consequence, we have that: 
	\[\|D^{\beta_\perp}\mathcal{E}^{k_s}_\rho\|_{L^2_{0,\eta}(\overline{\bf C}_0^\rho)} \le 
	C\rho^{\eta-\frac{1}{2}-k_\perp-k_s}\|\sigma\|_{H^{k_s}(\Lambda)}.\]
\end{lemma} 
\begin{proof}
	The first estimate follows directly from \eqref{estimate phi}, and from the application to $\sigma^{(\ell)}$ of the embedding $H^1(\Lambda)\subset L^{\infty}(\Lambda)$:
	\begin{align*}
		|D^{\beta_\perp}\mathcal{E}^{k_s}_\rho(r,\theta,s)| &= \bigg|D^{\beta_\perp}\phi_{2,\rho}(r)\sum_{\ell=0}^{k_s-1}\sigma^{(\ell)}(\rho)\phi^{k_s-1-\ell}_{1,\rho}(s-\rho) \,\dt\bigg| \\
		&\le C\rho^{-2-|\beta_\perp|} \sum_{\ell=0}^{k_s-1}\|\sigma^{(\ell)}\|_{H^1(\Lambda)}\rho^{\ell-k_s} \\
		&\le C\rho^{-2-k_\perp-k_s} \|\sigma\|_{H^{k_s}(\Lambda)}.
	\end{align*}
	For the second estimate, we apply the first one and then integrate in spherical coordinates as in the second estimate of the previous Lemma. We leave the details to the reader.  
\end{proof}

In order to treat the singularity at the extreme points of $\Lambda$ we will sometimes apply Lemma \ref{lemma Apqalpha} but with weight $d_e^{-\eta}$ instead of $d^{-\eta}$. The following remark is an analogous to Remark \ref{remark eta line} for this case. 

\begin{rmk}\label{rmk eta point}
	We consider Lemma \ref{lemma Apqalpha} with $E=\{\lambda(0)\}$ and weights $d_e^{-\eta}$ and $d_e^{-\eta^*}$. Combining \eqref{Apqalphacondeta} and \eqref{Apqalphacondeta*} with \eqref{eta*} and taking into account that $\alpha=2$ and  $\dim_A(\{\lambda(0)\})=0$ it is easy to check that the restrictions are reduced to: 
	\begin{equation}\label{cond eta point}
	\frac{1}{2}<\eta<\frac{3}{2},
	\end{equation}
	which gives a non-empty interval for $\eta$, for every $1<p\le q=2$.  
\end{rmk}

Finally we are now in possession of all the elements necessary to prove Theorem \ref{main theorem anisotropic}. We begin with the norm of $u_\rho^\circ$ in $C(\Lambda,\rho)$: 

\begin{prop}\label{anisotropic urho near}
	Given $\beta= \beta_s+\beta_\perp$ a multiindex, with $|\beta_s|=k_s$ and $|\beta_\perp|=k_\perp$. If conditions \eqref{gamma cond anisotropic} and \eqref{mu cond anisotropic} are fulfilled the following estimate holds 
	\[\|D^\beta u_\rho^\circ\|_{L^2_{\gamma,\mu}(C(\Lambda,2\rho))} \le C\|\sigma\|_{H^{k_s}(\Lambda)},\]
	with a constant $C$ independent of $\rho$. 
\end{prop}

Since the proof of this result is rather long, we split it in two lemmas. First, observe that
\begin{align*}
	\|D^{\beta_s+\beta_\perp} u_\rho^\circ\|_{L^2_{\gamma,\mu}(C(\Lambda,2\rho))} &\sim  \|I_2(D^{\beta} \sigma_\rho)\|_{L^2_{\gamma,\mu}(C(\Lambda,2\rho))},
\end{align*}
Recalling Lemma \ref{anisotropic deriv} we have that 
\[D^\beta \sigma_\rho = D^{\beta_\perp}(\sigma^{(k_s)})_\rho + D^{\beta_\perp}\mathcal{E}_\rho^{k_s}+D^{\beta_\perp}\mathcal{E}_{L-\rho}^{k_s},\]
and using this we get
	\begin{align*}
		\|(I_2(D^{\beta} \sigma_\rho)\|_{L^2_{\gamma,\mu}(C(\Lambda,2\rho))}   
		&\le \|I_2(D^{\beta_\perp} (\sigma^{(k_s)})_ \rho)\|_{L^2_{\gamma,\mu}(C(\Lambda,2\rho))} + \|I_2(D^{\beta_\perp} \mathcal{E}_\rho)\|_{L^2_{\gamma,\mu}(C(\Lambda,2\rho))}  \\ 
		&\quad\quad\quad\quad\quad\quad\quad\quad +  
		\|I_2(D^{\beta_\perp} \mathcal{E}^{k_s}_{L-\rho})\|_{L^2_{\gamma,\mu}(C(\Lambda,2\rho))}
	\end{align*}
It is clear that the second and third terms are completely analogous, so we devote the following lemmas to the estimation of the first and second terms.

\begin{lemma}\label{anisotropic sigma ks near}
	Under the conditions of Proposition \ref{anisotropic urho near} the following estimate holds
	\[\|I_2(D^{\beta_\perp} (\sigma^{(k_s)})_\rho)\|_{L^2_{\gamma,\mu}(C(\Lambda,2\rho))} \le C\|\sigma\|_{H^{k_s}(\Lambda)},\]
	where the constant $C$ is independent of $\rho$. 
\end{lemma}
\begin{proof}   
	If $\mu\ge 0$, we have that $d_e^{2\mu}(x)\le CL^{2\mu}$ for every $x\in C(\Lambda,2\rho)$. Hence, we can drop the weight $d_e$: 
	\[\|I_2(D^{\beta_\perp} (\sigma^{(k_s)})_\rho)\|_{L^2_{\gamma,\mu}(C(\Lambda,2\rho))}\le C \|I_2(D^{\beta_\perp} (\sigma^{(k_s)})_\rho)\|_{L^2_{\gamma}(C(\Lambda,2\rho))},\]
	and we can apply Lemma \ref{u_rho near} but with $(\sigma^{(k_s)})_\rho$ playing the role of $\sigma_\rho$ and $\beta_\perp$ playing the role of $\beta$, obtaining, under the assumption $\gamma>k_\perp-1$, 
	 \begin{align*}
		\|I_2(D^{\beta_\perp} (\sigma^{(k_s)})_\rho)\|_{L^2_{\gamma}(C(\Lambda,2\rho))} & \le C\|\sigma^{(k_s)}\|_{L^2(\Lambda)} \le C\|\sigma\|_{H^{k_s}(\Lambda)}.
	 \end{align*}
	 
	 The case $-\frac{1}{2}<\mu<0$ requires some additional effort, in order to handle the negative exponent. It is important to observe that this only occurs when $k_s=0$. Without loss of generality we assume that: 
	\[\|I_2(D^{\beta_\perp}\sigma_\rho)\|_{L^2_{\gamma,\mu}(C(\Lambda[\frac{L}{2},L],2\rho))} \le \|I_2(D^{\beta_\perp}\sigma_\rho)\|_{L^2_{\gamma,\mu}(C(\Lambda[0,\frac{L}{2}],2\rho))},\]
	so it is enough to estimate the norm in $C(\Lambda[0,\frac{L}{2}],2\rho)$. 

	Furthermore, we separate the norm, distinguishing the part that is near $\lambda(0)$ and the part that is far from it: 
	\begin{align*}
		\|I_2(D^{\beta_\perp}\sigma_\rho)\|_{L^2_{\gamma,\mu}(C(\Lambda[0,\frac{L}{2}],2\rho))}  
		&\le \|I_2(D^{\beta_\perp}\sigma_\rho)\|_{L^2_{\gamma,\mu}({\bf C}_0^{2\rho})} + \|I_2(D^{\beta_\perp}\sigma_\rho)\|_{L^2_{\gamma,\mu}(\cup_{j=1}^{J-1}{\bf C}_j^{2\rho})}=:I+II. 
	\end{align*}

	Let us begin considering $I$:
	\[I \le \|I_2(D^{\beta_\perp}\sigma_\rho \chi_{\overline{\bf C}_0^\rho})\|_{L^2_{\gamma,\mu}({\bf C}_0^{2\rho})} 
	+ \|I_2(D^{\beta_\perp}\sigma_\rho\chi_{\cup_{j=2}^J{\bf C}_j^\rho})\|_{L^2_{\gamma,\mu}({\bf C}_0^{2\rho})}=:I_A+I_B.\]
  Now, for the first of these terms, we can apply the dual characterization of the norm, and the H\"older inequality with weight, which gives:
	\begin{align*}
		I_A &= \sup_{g:\|g\|_{L^2_{-\gamma,-\mu}({\bf C}^{2\rho}_0)}=1} \int_{{\bf C}^{2\rho}_0} I_2(D^{\beta_\perp}\sigma_\rho \chi_{\overline{\bf C}^\rho_0})(x)g(x)\,\dx \\
			&= \sup_{g:\|g\|_{L^2_{-\gamma,-\mu}({\bf C}^{2\rho}_0)}=1}\int_{\overline{\bf C}^\rho_0} D^{\beta_\perp}\sigma_\rho(y) I_2(g\chi_{{\bf C}^{2\rho}_0})(y)\,\dy \\
		&\le \sup_{g:\|g\|_{L^2_{-\gamma,-\mu}({\bf C}^{2\rho}_0)}=1} \|D^{\beta_\perp}\sigma_\rho\|_{L^2_{0,\eta}(\overline{\bf C}^\rho_0)} \|I_2(g\chi_{{\bf C}^{2\rho}_0})\|_{L^2_{0,-\eta}(\overline{\bf C}^\rho_0)}. 
	\end{align*}
	For the first factor, we apply Lemma \ref{norm sigma_rho anisotropic}. 
 To the second factor we apply Lemma \ref{lemma Apqalpha} but with $E=\{\lambda(0)\}$ and weight $d_e^{-\eta}$. In particular, we choose, for some small value $\varepsilon>0$: 
 \[\eta = \frac{1}{2}+\varepsilon,\]
 which in turn implies (thanks to \eqref{eta*}): 
 \[\eta^*=\varepsilon+\frac{3}{p}-3.\]
 Then we multiply and divide by $d^{-\gamma}d_e^{-\mu}$ and apply H\"older's inequality with exponents $\frac{2}{p}$ and $\frac{2}{2-p}$:
	\begin{align*}
		\|I_2(g\chi_{{\bf C}_0^{2\rho}})\|_{L^2_{0,-\eta}(\overline{\bf C}_0^\rho)} &\le \|g\|_{L^p_{0,-\eta^*}({\bf C}_0^{2\rho})} \\
		&\le \|g\|_{L^2_{-\gamma,-\mu}({\bf C}_0^{2\rho})} \bigg(\int_{{\bf C}_0^{2\rho}}d(x)^{\gamma\frac{2p}{2-p}}d_e(x)^{(\mu-\eta^*)\frac{2p}{2-p}}\,\dx\bigg)^\frac{2-p}{2p}
	\end{align*}
	The norm of $g$ equals $1$. For the integral of the weights we enlarge the domain of integration to the ball $B(\lambda(0),\sqrt{20}\rho)$, and integrate using spherical coordinates $(r,\theta,\xi)$ recalling Remark \ref{two weight trick}: 
	\begin{align*}
		\|I_2(g\chi_{{\bf C}_0^{2\rho}})\|_{L^2_{0,-\eta}(\overline{\bf C}_0^\rho)}  &\le C\bigg(\int_{{\bf C}_0^{2\rho}}d(x)^{\gamma\frac{2p}{2-p}}d_e(x)^{(\mu-\eta^*)\frac{2p}{2-p}}\,\dx\bigg)^\frac{2-p}{2p}\\
		&\le C \bigg(\int_{0}^{\sqrt{20}\rho}\int_0^{2\pi}\int_0^{\frac{\pi}{2}}
		r^{(\gamma+\mu-\eta^*)\frac{2p}{2-p}+2}\sin(\xi)^{\gamma\frac{2p}{2-p}+1}
		\,\dxi\dtheta\dr\bigg)^\frac{2-p}{2p} \\
		&\le C \bigg(\int_{0}^{\sqrt{20}\rho}
		r^{(\gamma+\mu-\eta^*)\frac{2p}{2-p}+2}\dr
		\,\int_0^{\frac{\pi}{2}}\sin(\xi)^{\gamma\frac{2p}{2-p}+1}\dxi\bigg)^\frac{2-p}{2p} \\
		&\le C \rho^{\gamma+\mu-\eta^*+3\frac{2-p}{2p}},
	\end{align*}
	where in the last step we assumed integrability conditions on both integrals. For the first one, we need $(\mu+\gamma-\eta^*)\frac{2p}{2-p}+2>-1$ which is equivalent to: 
	\[\gamma+\mu>\frac{3}{2}-\frac{3}{p}+\eta^* = -\frac{3}{2}+\varepsilon.\] 
	Since $\mu+\gamma>-\frac{3}{2}$, we can choose $\varepsilon$ smal enough such that the integrability condition in fulfilled. For the second integral, we need to impose the condition $\gamma\frac{2p}{2-p}+1>-1$, which is equivalent to 
	\[\gamma>1-\frac{2}{p}.\] 
	And now, since $\gamma>-1$, we choose $p$ close enough to $1$ so that the condition is satisfied. 

	Finally, we join both factors and apply \eqref{gamma cond anisotropic} and \eqref{mu cond anisotropic} obtaining the desired estimate. 
	\begin{align*}
		I_A &\le C\rho^{\eta-\frac{1}{2}-k_\perp}\rho^{\gamma+\mu-\eta^*+3\frac{2-p}{2p}}\|\sigma\|_{H^1(\Lambda)} \\
			&\le  C\rho^{\varepsilon-k_\perp +\gamma+\mu-\varepsilon-\frac{3}{p}+3+\frac{3}{p}-\frac{3}{2}}\|\sigma\|_{H^1(\Lambda)} \\
		&\le C\rho^{\gamma + 1 -k_\perp + \mu +\frac{1}{2}}\|\sigma\|_{H^1(\Lambda)}\\
		&\le C\|\sigma\|_{H^1(\Lambda)}.
	\end{align*}

	For $I_B$, we use that for  $y\in {\bf C}_j^{\rho}$ with $j\ge 2$ and $x\in {\bf C}_0^{2\rho}$, $|x-y|\sim 2^j\rho$, so applying the first estimate in Lemma \ref{norm sigma_rho anisotropic} we obtain 
	\begin{align*}
		\sum_{j=2}^{J} \int_{{\bf C}_j^\rho} \frac{1}{|x-y|}D^{\beta_\perp}\sigma_\rho(y)\,\dy &\le \sum_{j=2}^{J}(2^j\rho)^{-1} \int_{{\bf C}_j^\rho} D^{\beta_\perp}\sigma_\rho(y)\,\dy \\				
	  &\le C \sum_{j=2}^{J}(2^j\rho)^{-1} \rho^{-2-k_\perp}\|\sigma\|_{H^1(\Lambda)}|{\bf C}_j^\rho| \\ 
		&\le C\rho^{-k_\perp}J \|\sigma\|_{H^1(\Lambda)}. 
	\end{align*}
	Inserting this in $I_B$ and recalling that $J\sim |\log(\rho)|$, we obtain 
	\begin{align*}
		I_B^2 &= \int_{{\bf C}_0^{2\rho}} \Big|\sum_{j=2}^J\int_{{\bf C}_j^{2\rho}}\frac{1}{|x-y|}D^{\beta_\perp}\sigma_\rho(y)\,\dy\Big|^2 d(x)^{2\gamma}d_e(x)^{2\mu}\,\dx \\
					&\le C\rho^{-2k_\perp}\log^2(\rho)\|\sigma\|_{H^1(\Lambda)}^2 \int_{{\bf C}_0^{2\rho}} d(x)^{2\gamma}d_e(x)^{2\mu}\,\dx
	\end{align*}
	For the integral of the weights we apply the argument of Remark \ref{two weight trick}. Taking spherical coordinates $(r,\theta,\xi)$ on a ball containing ${\bf C}_0^{2\rho}$ and using the integrability conditions \eqref{gamma cond anisotropic} and \eqref{mu cond anisotropic} we get 
	\begin{align*}
		I_B &\le C\rho^{-k_\perp}|\log(\rho)|\|\sigma\|_{H^1(\Lambda)}\bigg(\int_0^{\sqrt{5}\rho} \int_0^{\frac{\pi}{2}}\int_0^{2\pi}  r^{2\mu+2\gamma+2}\sin(\xi)^{2\gamma+1}\,\dtheta\dxi\dr\bigg)^\frac{1}{2} \\
				&\le C\rho^{-k_\perp}|\log(\rho)| \|\sigma\|_{H^1(\Lambda)}\bigg(\int_0^{\sqrt{5}\rho} r^{2\mu+2\gamma+2}\,\dr \int_0^{\frac{\pi}{2}} \sin(\xi)^{2\gamma+1}\,\dxi\bigg)^\frac{1}{2} \\
				&\le C\rho^{\mu+\gamma+\frac{3}{2}-k_\perp}|\log(\rho)| \|\sigma\|_{H^1(\Lambda)}^2,
	\end{align*}
	 Moreover, the exponent of $\rho$ is positive thanks to \eqref{gamma cond anisotropic} and \eqref{mu cond anisotropic} and consequently we have that 
	\[\rho^{\mu+\gamma+\frac{3}{2}-k_\perp}|\log(\rho)|\le C,\]
	for $\rho$ approaching $0$, which concludes the estimate for $I_2$. 

	Finally, we consider the term $II$. Taking into account that for $x\in {\bf C}_j^{2\rho}$ with $j\ge 1$, $d_e(x)\sim 2^j\rho$ we have that 
	\begin{align*}
		II^2 &=\sum_{j=1}^{J-1}\int_{{\bf C}_j^{2\rho}}|I_2(D^\beta \sigma_\rho)(x)|^2d(x)^{2\gamma}d_e(x)^{2\mu}\dx \\
				 &\le C \sum_{j=1}^{J-1}(2^j\rho)^{2\mu} \int_{{\bf C}_j^{2\rho}}|I_2(D^\beta \sigma_\rho)(x)|^2d(x)^{2\gamma}\dx. 
	\end{align*}

	Now we estimate the integral inside the summation, which is the squared $L^2_\gamma$ norm of $I_2(D^\beta \sigma_\rho)$ over ${\bf C}_j^{2\rho}$. We separate this norm into three parts, localizing $\sigma_\rho$ in a neighbourhood of ${\bf C}_j^{2\rho}$, and far from it: 
	\begin{align*}
		\|I_2(D^\beta \sigma_\rho)\|_{L^2_\gamma({\bf C}_j^{2\rho})} &\le \|I_2(D^\beta \sigma_\rho \chi_{\overline{\bf C}_j^{2\rho}})\|_{L^2_\gamma({\bf C}_j^{2\rho})} \\ 
			&\quad\quad\quad+ \|I_2(D^\beta \sigma_\rho\chi_{\cup_{i=0}^{j-2}{\bf C}_i^{2\rho}})\|_{L^2_\gamma({\bf C}_j^{2\rho})} + \|I_2(D^\beta \sigma_\rho\chi_{\cup_{i=j+2}^J}{\bf C}_i^{2\rho})\|_{L^2_\gamma({\bf C}_j^{2\rho})}\\
			&=: II_A + II_B + II_C.
		\end{align*}
		Naturally, the second term vanishes if $j=1$ whereas the third one vanishes if $j=J-1$. 

		In $II_A$ we apply Lemma \ref{lemma Apqalpha} with $E=\Lambda$. We begin, as usual, by using the dual characterization of the norm, applying Fubini's lemma and the Cauchy-Schwartz inequality: 
		\begin{align*}
			II_A &= \sup_{g:\|g\|_{L^2_{-\gamma}({\bf C}_j^{2\rho})}=1} \int_{{\bf C}_j^{2\rho}}I_2(D^{\beta_\perp}\sigma_\rho\chi_{\overline{\bf C}_j^\rho})(x)g(x)\dx \\
					 &=\sup_{g:\|g\|_{L^2_{-\gamma}({\bf C}_j^{2\rho})}=1} \int_{\overline{\bf C}_j^{\rho}}D^{\beta_\perp}\sigma_\rho(y) I_2(g\chi_{{\bf C}_j^{2\rho}})(y)\dx \\ 
					 &=\sup_{g:\|g\|_{L^2_{-\gamma}({\bf C}_j^{2\rho})}=1} \|D^{\beta_\perp}\sigma_\rho\|_{L^2_\eta(\overline{\bf C}_j^{\rho})}  \|I_2(g\chi_{{\bf C}_j^{2\rho}})\|_{L^2_{-\eta}(\overline{\bf C}_j^{\rho})} \\
					 &\le C \sup_{g:\|g\|_{L^2_{-\gamma}({\bf C}_j^{2\rho})}=1} \|D^{\beta_\perp}\sigma_\rho\|_{L^2_\eta(\overline{\bf C}_j^{\rho})}  \|g\|_{L^p_{-\eta^*}({\bf C}_j^{2\rho})}
		\end{align*}
		The first factor is estimated by Lemma \ref{norm sigma_rho anisotropic}. For the second one we choose $\eta$ and $p$ as in Lemma \ref{u_rho near}. Then, we multiply and divide by $d^\gamma$, we apply the H\"older inequality with exponents $2/p$ and $2/(2-p)$, and finally we integrate the weight in cylindrical coordinates 
		\begin{align*}
			\|g\|_{L^p_{-\eta^*}({\bf C}_j^{2\rho})} &\le \|g\|_{L^2_{-\gamma}({\bf C}_j^{2\rho})} \bigg(\int_{{\bf C}_j^{2\rho}}d(x)^{(\gamma-\eta^*)\frac{2p}{2-p}}\dx\bigg)^\frac{2-p}{2p} \\
				&\le C(2^j\rho)^{\frac{2-p}{2p}} \bigg(\int_0^{2\rho} r^{(\gamma-\eta^*)\frac{2p}{2-p}+1}\dr\bigg)^\frac{2-p}{2p} \\
				&\le C (2^j\rho)^{\frac{1}{p}-\frac{1}{2}}\rho^{\gamma-\eta^*+\frac{2-p}{p}}.
		\end{align*}
		Joining both estimates and recalling from Lemma \ref{u_rho near} that $\eta-\eta^*=\frac{7}{2}-\frac{3}{p}$ we obtain: 
	\begin{align*}
		II_A &\le C (2^j\rho)^\frac{1}{p} \rho^{\gamma+\eta-\eta^*-1-k_\perp+\frac{2}{p}-1}\|\sigma\|_{H^1(\Lambda)} \\
				 &\le (2^j\rho)^\frac{1}{p} C \rho^{\gamma+\frac{7}{2}-\frac{3}{p}-2-k_\perp +\frac{2}{p}}\|\sigma\|_{H^1(\Lambda)}\\
				 &\le C(2^j\rho)^\frac{1}{p}  \rho^{\gamma +\frac{3}{2}-k_\perp-\frac{1}{p}}\|\sigma\|_{H^1(\Lambda)}.
	\end{align*}

	For $II_B$, let us observe that if $x\in{\bf C}_j^{2\rho}$ and $y\in {\bf C}_i^\rho$ with $i<j-1$, then: $|x-y|\sim 2^j\rho-2^i\rho \sim 2^j\rho$. Consequently, applying the first estimate in Lemma \ref{norm sigma_rho anisotropic} we have
\begin{align*}
|I_2(D^{\beta_\perp}\sigma_\rho \chi_{\cup_{i=0}^{j-1}{\bf C}_i^\rho})(x)| &= \bigg|\sum_{i=0}^{j-1}\int_{{\bf C}_i^\rho}|x-y|^{-1} D^{\beta_\perp}\sigma_\rho(y)\dy\bigg| \\
		&\le C(2^j\rho)^{-1}  \sum_{i=0}^{j-1} \int_{{\bf C}_i^\rho} |D^{\beta_\perp}\sigma(y)|\dy\\
		&\le C(2^j\rho)^{-1} \rho^{-2-k_\perp}\|\sigma\|_{H^1(\Lambda)}\sum_{i=0}^{j-1} (2^i\rho)\\
		&\le C(2^j\rho)^{-1}(2^j\rho) \rho^{-2-k_\perp}\|\sigma\|_{H^1(\Lambda)} \\
&\le C\rho^{-2-k_\perp}\|\sigma\|_{H^1(\Lambda)}
	\end{align*}
Inserting this estimate in $II_B$ and integrating in cylindrical coordinates we obtain
\begin{align*}
	II_B &\le  C\rho^{-2-k_\perp}\|\sigma\|_{H^1(\Lambda)}\bigg(\int_{{\bf C}_j^{2\rho}}d(x)^{2\gamma}\dx\bigg)^\frac{1}{2} \\
			 &\le C\rho^{-2-k_\perp}\|\sigma\|_{H^1(\Lambda)} \rho^{\gamma+1}(2^j\rho)^\frac{1}{2} \\
			 &\le C(2^j\rho)^\frac{1}{2}\rho^{\gamma-1-k_\perp}\|\sigma\|_{H^1(\Lambda)}
\end{align*}

The argument for $II_C$ is essentially the same, but taking into account that if $x\in{\bf C}_j^{2\rho}$ and $y\in {\bf C}_i^\rho$ with $i>j+1$, then: $|x-y|\sim 2^i\rho-2^j\rho \sim 2^i\rho$. Following the analysis carried our for $II_B$, and recalling that $J\sim |\log(\rho)|$, this leads us to: 
\begin{align*}
	|I_2(D^{\beta_\perp}\sigma_\rho \chi_{\cup_{i=0}^{j-1}{\bf C}_i^\rho})(x)| &\le C \rho^{-2-k_\perp}\|\sigma\|_{H^1(\Lambda)}\sum_{i=j+1}^{J} (2^i\rho)^{-1}(2^i\rho) \\
  &\le C \rho^{-2-k_\perp}\|\sigma\|_{H^1(\Lambda)}J \\
  &\le C \rho^{-2-k_\perp}|\log(\rho)|\|\sigma\|_{H^1(\Lambda)},
\end{align*}
and inserting this in $II_C$ and integrating in cylindrical coordinates we obtain
\begin{align*}
 II_C &\le C \rho^{-2-k_\perp} |\log(\rho)|\|\sigma\|_{H^1(\Lambda)} \bigg(\int_{{\bf C}_j^{2\rho}}d(x)^{2\gamma}\dx\bigg)^\frac{1}{2} \\
			&\le C(2^j\rho)^\frac{1}{2} \rho^{\gamma-1-k_\perp} |\log(\rho)| \|\sigma\|_{H^1(\Lambda)}.
\end{align*}
Since $p<2$ we can estimate: 
\[II_A+II_B+II_C \le C (2^j\rho)^\frac{1}{p}\rho^{\gamma+\frac{3}{2}-k_\perp-\frac{1}{p}}|\log(\rho)| \|\sigma\|_{H^1(\Lambda)}.\]
And with this we can finally complete the estimate for $II$:
\begin{align*}
	II &\le C \bigg(\sum_{j=1}^{J-1} (2^j\rho)^{2\mu}(II_A+II_B+II_C)^2\bigg)^\frac{1}{2}\\
		 &\le C\rho^{\gamma+\frac{3}{2}-k_\perp-\frac{1}{p}}|\log(\rho)| \|\sigma\|_{H^1(\Lambda)}\bigg(\sum_{j=1}^{J-1} (2^j\rho)^{2\mu+\frac{2}{p}}\bigg)^\frac{1}{2}.
\end{align*}
The summation is bounded by $(2^J\rho)^{2\mu+\frac{2}{p}}\sim L^{2\mu+\frac{2}{p}}$, which is a constant independent of $\rho$. Moreover, since $\gamma>k_\perp-1$, $p$ can be chosen as close to $2$ as needed such that $\gamma+\frac{3}{2}-k_\perp-\frac{1}{p}> 0$, and consequently the factor $\rho^{\gamma+\frac{3}{2}-k_\perp-\frac{1}{p}}|\log(\rho)|$ is bounded for $\rho$ tending to $0$. Hence,
\[II\le C \|\sigma\|_{H^1(\Lambda)}.\]
This completes the proof. 
	\end{proof}

	\begin{lemma}\label{anisotropic E near}
	Under the conditions of Proposition \ref{anisotropic urho near} the following estimate holds
 \[\|I_2(D^{\beta_\perp} \mathcal{E}^{k_s}_\rho)\|_{L^2_{\gamma,\mu}(C(\Lambda,2\rho))} \le C\|\sigma\|_{H^{k_s}(\Lambda)},\]
 where the constant $C$ is independent of $\rho$. 
\end{lemma}
\begin{proof}
	We only consider the norm of $I_2(D^{\beta_\perp} \mathcal{E}_\rho)$ over $C(\Lambda[0,\frac{L}{2}],2\rho)$. It is clear that the norm over $C(\Lambda[\frac{L}{2},L],2\rho)$ can be estimated by means of the same arguments. We denote $\widetilde{\bf C}_0^\delta = \overline{\bf C}_0^\delta\cup {\bf C}_2^\delta$. Then:
	 \begin{align*}
		\|I_2(D^{\beta_\perp} &\mathcal{E}^{k_s}_\rho)\|_{L^2_{\gamma,\mu}(C(\Lambda[0,\frac{L}{2}],2\rho))}\\ 
		&\le \|I_2(D^{\beta_\perp} \mathcal{E}^{k_s}_\rho)\|_{L^2_{\gamma,\mu}(\widetilde{\bf C}_0^{2\rho})} + \|I_2(D^{\beta_\perp} \mathcal{E}^{k_s}_\rho)\|_{L^2_{\gamma,\mu}(\cup_{j=3}^{J-1}{\bf C}_j^{2\rho})}\\
		&=: I + II.
	 \end{align*}
	 
	 For $I$ we use the dual characterization of the norm and apply Fubini's Lemma and Cauchy-Schwartz's inequality obtaining:
		\begin{align*}I &= \sup_{g:\|g\|_{L^2_{-\gamma,-\mu}(\widetilde{\bf C}_0^{2\rho})}=1} \int_{\widetilde{\bf C}_0^{2\rho}}  I_2(D^{\beta_\perp} \mathcal{E}^{k_s}_\rho)(x) g(x)\,\dx \\
		&\le \sup_{g:\|g\|_{L^2_{-\gamma,-\mu}(\widetilde{\bf C}_0^{2\rho})}=1} \int_{\overline{\bf C}_0^{2\rho}} D^{\beta_\perp} \mathcal{E}^{k_s}_\rho(y) I_2(g\chi_{\widetilde{\bf C}_0^{2\rho}})(y) \,\dy \\
		&\le \sup_{g:\|g\|_{L^2_{-\gamma,-\mu}(\widetilde{\bf C}_0^{2\rho})}=1} \| D^{\beta_\perp} \mathcal{E}^{k_s}_\rho\|_{L^2_{0,\eta}(\overline{\bf C}_0^{2\rho})}\|I_2(g\chi_{\widetilde{\bf C}_0^{2\rho}})\|_{L^2_{0,-\eta}(\overline{\bf C}_0^{2\rho})}
		\end{align*}
		
		The first factor in the supremum is bounded by Lemma \ref{E estimate}, whereas for the second one we proceed exactly as in the estimate of the term $I_A$ in the previous lemma, obtaining: 

		\[\|I_\alpha(g\chi_{\widetilde{\bf C}_0^{2\rho}})\|_{L^2_{0,-\eta}(\overline{\bf C}_0^{2\rho})} \le C\rho^{\gamma+\mu-\eta^*+3\frac{2-p}{2p}}.\]

		Joining both estimates, and using as in the estimate of $I_A$ in the previous lemma that $\eta-\eta^*=\frac{7}{2}-\frac{1}{p}$ we have: 
		\begin{align*}
			I &\le C\rho^{\eta-\frac{1}{2}-k_s-k_\perp+\gamma+\mu-\eta^*+3\frac{2-p}{2p}}\|\sigma\|_{H^{k_s}(\Lambda)} \\
			&\le C\rho^{\gamma+\mu+\frac{3}{2}-k_s-k_\perp}\|\sigma\|_{H^{k_s}(\Lambda)},
		\end{align*}
		and the exponent of $\rho$ is positive for $\gamma>k_\perp-1$ and $\mu>k_s-\frac{1}{2}$.

		For $II$ we begin by integrating by parts passing the derivatives with respect to $s$ from $\phi_{1,\rho}$ to the kernel $\Gamma$. We denote $\partial_{s_y}^{k_s-1-\ell}$ the derivative with respect to $s_y$ of order $k_s-1-\ell$:
\begin{align*}
	|I_2(D^{\beta_\perp}\mathcal{E}^{k_s}_\rho)| &= \Big|\sum_{\ell=0}^{k_s-1}\int_{\overline{\bf C}_0^\rho} \Gamma(x-y) \sigma^{(\ell)}(\rho) \phi_{1,\rho}^{(k_s-1-\ell)}(s_y-\rho) D^{\beta''_\perp}\phi_{2,\rho}(r_y) \,\dy\Big| \\ 
	&\le \sum_{\ell=0}^{k_s-1}\int_{\overline{\bf C}_0^\rho} \Big| \partial_{s_y}^{k_s-1-\ell} \Gamma(x-y) \sigma^{(\ell)}(\rho) \phi_{1,\rho}(s_y-\rho) D^{\beta_\perp}\phi_{2,\rho}(r_y)\Big| \,\dy \\ 
	&\le \sum_{\ell=0}^{k_s-1}\int_{\overline{\bf C}_0^\rho}  |x-y|^{\ell-k_s} |\sigma^{(\ell)}(\rho)| \phi_{1,\rho}(s_y-\rho) |D^{\beta_\perp}\phi_{2,\rho}(r_y)| \,\dy,
\end{align*}
where we used that $\phi_{1,\rho}(s-\rho)\phi_{2,\rho}(r)$ and its derivatives vanish at the boundary of $\overline{\bf C}_0^\rho$, and consequently so do the boundary terms from the integration by parts. We now apply the estimates \eqref{estimate phi} and \eqref{embedding}, which give
\begin{equation}
	|I_2(D^{\beta_\perp}\mathcal{E}^{k_s}_\rho)| \le C \rho^{-3-k_\perp}\|\sigma\|_{H^{k_s}(\Lambda)}\int_{\overline{\bf C}_0^\rho}  |x-y|^{-k_s} \,\dy,
\end{equation}

Now, we insert this estimate in the norm. For every $x\in {\bf C}_j^{2\rho}$ with $j\ge 3$, $d_e(x)\sim 2^j\rho$. Moreover, for $y\in {\bf C}_0^\rho$, $|x-y|\sim 2^j\rho$. Using these and integrating in cylindrical coordinates we obtain: 
\begin{align*}
	II &= \|I_2(D^{\beta_\perp} \mathcal{E}^{k_s}_\rho)\|_{L^2_{\gamma,\mu}(\widetilde{\bf C}_0^{2\rho})} \\
 &\le C \rho^{-3-k_\perp}\|\sigma\|_{H^{k_s}(\Lambda)} \bigg\{\sum_{j=3}^{J-1} \int_{{\bf C}_j^{2\rho}} \Big|\int_{\overline{\bf C}_0^\rho}  |x-y|^{-k_s} \,\dy\Big|^2 d(x)^{2\gamma}d_e(x)^{2\mu}\,\dx\bigg\}^\frac{1}{2} \\
 &\le C \rho^{-3-k_\perp}\|\sigma\|_{H^{k_s}(\Lambda)} \bigg\{\sum_{j=3}^{J-1} (2^j\rho)^{2\mu-2k_s} |\overline{\bf C}_0^{\rho}|^2\int_{{\bf C}_j^{2\rho}}  d(x)^{2\gamma}\,\dx\bigg\}^\frac{1}{2} \\
	&\le C \rho^{-3-k_\perp}\rho^3\|\sigma\|_{H^{k_s}(\Lambda)} \bigg\{\sum_{j=3}^{J-1} (2^j\rho)^{2\mu-2k_s} \int_{2^j\rho}^{2^{j+1}\rho} \int_0^{2\pi}\int_0^{2\rho}r^{2\gamma+1}\,\dr\dtheta\ds\bigg\}^\frac{1}{2} \\
	&\le C \rho^{-k_\perp}\|\sigma\|_{H^{k_s}(\Lambda)} \bigg\{\sum_{j=2}^{J-1} (2^j\rho)^{2\mu-2k_s} (2^j\rho) \rho^{2\gamma+2}\bigg\}^\frac{1}{2} \\  
	&\le C \rho^{\gamma+1-k_\perp}\|\sigma\|_{H^{k_s}(\Lambda)} \bigg\{\sum_{j=2}^{J-1} (2^j\rho)^{2\mu-2k_s+1}\bigg\}^\frac{1}{2}.
\end{align*}

Since $\mu>k_s-\frac{1}{2}$, $2\mu-2k_s+1>0$ so the summation on $j$ is $\sim (2^J\rho)^{2\mu-2k_s+1} \sim L^{2\mu-2k_s+1}$. Likewise, since $\gamma>k_\perp-1$, $\rho^{\gamma+1-k_\perp}\le 1$, which yields
\[II\le C\|\sigma\|_{H^k(\Lambda)},\]
completing the proof. 
\end{proof}

The previous two lemmas constitute the proof of Proposition \ref{anisotropic urho near}. Let us now complete the analysis of $u_\rho^\circ$ in a close neighbourhood of $\Lambda$ by estimating its norm over $B^+_\rho(0)\cup B^+_\rho(L)$. 
\begin{lemma}\label{anisotropic urho near ball}
	Let $\beta$ be a multiindex and $k=|\beta|$. If $\mu>k-\frac{3}{2}$, then:
	\[\|D^\beta u_\rho^\circ\|_{L^2_{0,\mu}(B^+_{2\rho}(0)\cup B^+_{2\rho}(L))}\le C\|\sigma\|_{H^m(\Lambda)},\]
	with a constant $C$ independent of $\rho$.
\end{lemma}
\begin{proof}
	In $B^+_\rho(0)\cup B^+_\rho(L)$, $d(x)=d_e(x)$, so the weight is $d_e(x)^{2(\gamma+\mu)}$. 
Thanks to Lemma \ref{anisotropic deriv}, we have that: 
	\[\|D^\beta u_\rho^\circ\|_{L^2_{\gamma,\mu}(B^+_{2\rho}(0))} \le \|D^{\beta_\perp}(\sigma^{(k_s)})_\rho\|_{L^2_{\gamma,\mu}}(B^+_{2\rho}(0)) + \|D^{\beta_\perp}\mathcal{E}^{k_s}_\rho\|_{L^2_{\gamma,\mu}}(B^+_{2\rho}(0)).\]
	The estimate for the first term is almost exactly the same than the one given for the term $I$ in the proof of Lemma \ref{anisotropic sigma ks near}. The only minor difference lies in the fact that in that lemma we were integrating over ${\bf C}_\rho^0$, where both weights appear. Here we only need to consider the weight $d_e(x)^{2(\gamma+\mu)}$, for which the integrability condition is $\mu+\gamma>-\frac{3}{2}$. Moreover, following Lemma \ref{anisotropic sigma ks near}, at the end of the estimate we obtain a factor $\rho^{\gamma+\mu+\frac{3}{2}-k_\perp}|\log(\rho)|$ which is bounded for $\rho\to 0$ whenever $\mu+\gamma>k_\perp-\frac{3}{2}$. 

	For the second term the situation is quite similar. Indeed, the estimate is analogous to the one for the term $I$ in Lemma \ref{anisotropic E near}. Once again the weight $d(x)$ is absorbed by $d_e(x)$. At the end we have a factor $\rho^{\gamma+\mu+\frac{3}{2}-k}$ which is bouned for $\gamma+\mu>k-\frac{3}{2}$. 
\end{proof}

As in the previous section we now proceed to estimate the norms in $C(\Lambda,R_0)\setminus C(\Lambda,\rho)$. In particular, we prove: 

\begin{prop}\label{prop urho anisotropic not so near}
	Given $\beta=\beta_s+\beta_\perp$ a multi-index con $|\beta_s|=k_s$ and $|\beta_\perp|=k_\perp$. There is constant $C$ independent of $\rho$ such that 
	\[\|D^\beta u_\rho^\circ\|_{L^2_{\gamma,\mu}(C(\Lambda,R_0)\setminus C(\Lambda,2\rho))} \le C \|\sigma\|_{H^{k_s}(\Lambda)},\]
	for every $\gamma$ and $\mu$ satisfying \eqref{gamma cond anisotropic} and \eqref{mu cond anisotropic} respectively.   
\end{prop} 
 
The proof of this result is a combination of the arguments in Proposition \ref{anisotropic urho near} and the techniques of Lemma \ref{u_rho not so near}. As we did in Proposition \ref{anisotropic urho near}, we separate the norm into three parts:

\[\begin{split}
	\|D^\beta u_\rho^\circ&\|_{L^2_{\gamma,\mu}(C(\Lambda,R_0)\setminus C(\Lambda,2\rho))}\le \|D^{\beta_\perp}\Gamma*D^{\beta_s}\sigma_\rho\|_{L^2_{\gamma,\mu}(C(\Lambda,R_0)\setminus C(\Lambda,2\rho))} \\
												&\le \|(D^{\beta_\perp}\Gamma)*(\sigma^{(k_s)})_\rho\|_{L^2_{\gamma,\mu}(C(\Lambda,R_0)\setminus C(\Lambda,2\rho))} + \|(D^{\beta_\perp}\Gamma)*(\mathcal{E}_\rho^{k_s})\|_{L^2_{\gamma,\mu}(C(\Lambda,R_0)\setminus C(\Lambda,2\rho))}\\
												&\quad+ \|(D^{\beta_\perp}\Gamma)*(\mathcal{E}_{L-\rho}^{k_s})\|_{L^2_{\gamma,\mu}(C(\Lambda,R_0)\setminus C(\Lambda,2\rho))}
\end{split}\]

Once again, the second and third terms are completely analogous, so we devote the following lemmas to the study of the first and second terms. 

\begin{lemma}\label{anistropic sigma ks not so near}
	Under the conditions of Proposition \ref{anisotropic urho near}, there is a constant $C$ independent of $\rho$ such that
\[\|(D^{\beta_\perp}\Gamma)*(\sigma^{(k_s)})_\rho\|_{L^2_{\gamma,\mu}(C(\Lambda,R_0)\setminus C(\Lambda,2\rho))} \le C\|\sigma\|_{H^{k_s}(\Lambda)}.\]
\end{lemma}
\begin{proof}
	First, we observe that if $\mu\ge 0$ $d_e(x)^{2\mu}\le ((\frac{L}{2})^2+R_0^2)^{\mu}$, and we have that: 
	\[\|(D^{\beta_\perp}\Gamma)*(\sigma^{(k_s)})_\rho\|_{L^2_{\gamma,\mu}(C(\Lambda,R_0)\setminus C(\Lambda,2\rho))}\le \|(D^{\beta_\perp}\Gamma)*(\sigma^{(k_s)})_\rho\|_{L^2_{\gamma}(C(\Lambda,R_0)\setminus C(\Lambda,2\rho))},\]
	and the result follows from the application of Lemma \ref{u_rho not so near} with $\sigma^{(k_s)}$ in the place of $\sigma$. Hence, we only need to consider the case $-\frac{1}{2}<\mu<0$ which only occurs if $k_s=0$. In this case, we have that
	\begin{align*}
		|(D^{\beta_\perp}\Gamma)*(\sigma_\rho)(x)| &\le \int_{C(\Lambda,\rho)}\frac{1}{|x-y|^{1+k_\perp}}|\sigma(y)|\,\dy 
	\end{align*}
	Without loss of generality we may assume that
	\[\|(D^{\beta_\perp}\Gamma)*(\sigma_\rho)\|_{L^2_{\gamma,\mu}(C(\Lambda[\frac{L}{2},L],R_0)\setminus C(\Lambda,2\rho))}\le \|(D^{\beta_\perp}\Gamma)*(\sigma_\rho)\|_{L^2_{\gamma,\mu}(C(\Lambda[0,\frac{L}{2}],R_0)\setminus C(\Lambda,2\rho))},\]
	so we only need to estimate the norm over $C(\Lambda[0,\frac{L}{2}],R_0)\setminus C(\Lambda[0,\frac{L}{2}],2\rho)$. 

	Following the proof of Lemma \ref{u_rho not so near}, we have that
	\[C(\Lambda,\rho)=\{y\in C(\Lambda,\rho):|s_y-s_x|<r_x\}\cup\{y\in C(\Lambda,\rho):|s_y-s_x|\ge r_x\} =: A\cup B,\]
	and 
	\[\int_{C(\Lambda,\rho)}\frac{1}{|x-y|^{1+k_\perp}}|\sigma(y)|\,\dy = \int_{A}\frac{1}{|x-y|^{1+k_\perp}}|\sigma(y)|\,\dy + \int_{B}\frac{1}{|x-y|^{1+k_\perp}}|\sigma(y)|\,\dy =: I + II.\]
  Moreover,  
\[I\le C r_x^{-1-k_\perp}\int_{t:|t-s_x|<r_x+\rho}|\sigma(t)|\,\dt.\]
And now we can apply the compact embedding $H^1(\Lambda)\subset L^\infty(\Lambda)$, which gives: 
\[I\le C r_x^{-k_\perp}\|\sigma\|_{H^1(\Lambda)}.\]

Hence, since in the cylinder $C(\Lambda,R_0)$, $r_x=d(x)$, we have that
\begin{align*}
	\|I\|_{L^2_{\gamma,\mu}(C(\Lambda[0,\frac{L}{2}],R_0)\setminus C(\Lambda,2\rho))}^2 &\le C \|\sigma\|_{H^1(\Lambda)}^2 \int_{C(\Lambda[0,\frac{L}{2}],R_0)\setminus C(\Lambda,2\rho)} r_x^{-2k_\perp}d(x)^{2\gamma}d_e(x)^{2\mu}\,\dx \\
 &= C \|\sigma\|_{H^1(\Lambda)}^2 \int_{C(\Lambda[0,\frac{L}{2}],R_0)\setminus C(\Lambda,2\rho)} d(x)^{2\gamma-2k_\perp}d_e(x)^{2\mu}\,\dx.
\end{align*}
The integral can be estimated by enlarging the domain of integration to a ball containing $C(\Lambda,\frac{L}{2},R_0)$ and integrating in spherical coordinates, taking into account Remark \ref{two weight trick}:
\begin{align*}
	\int_{C(\Lambda[0,\frac{L}{2}],R_0)\setminus C(\Lambda,2\rho)} d(x)^{2\gamma-2k_\perp}d_e(x)^{2\mu}\,\dx
	&\le \int_0^{2\pi}\int_0^{\pi} \int_0^L r^{2\gamma+2\mu-2k_\perp+2}\sin(\xi)^{2\gamma-2k_\perp+1}\dr\dxi\dtheta \\
	&\le C \int_{0}^{\frac{\pi}{2}}\sin(\xi)^{2\gamma-2k_\perp+1}\dxi \int_{2\rho}^\frac{L}{2}r^{2\gamma+2\mu-2k_\perp+2}\dr \le C.
\end{align*}
For computing the integrals we have used the integrability conditions: $2\gamma-2k_\perp+1>-1$ and $2\gamma+2\mu-2k_\perp+2>-1$ which are satisfied thanks to \eqref{gamma cond anisotropic} and \eqref{mu cond anisotropic}. This completes the estimate for $I$.  

For $II$, once again we follow the estimate in the proof of Lemma \ref{u_rho not so near}. We begin with the decomposition: 

\[B = \bigcup_{j=0}^{J-1}\Big\{y\in C(\Lambda,\rho):\;2^jr_x<|s_y-s_x|<2^{j+1}r_x\Big\} =: \bigcup_{j=0}^{J-1}B_j,\]
where $J$ is the first integer such that $2^Jr_x>\frac{L}{2}$. 

Now, we estimate $II$ as in Lemma \ref{u_rho not so near}, obtaining: 
\[II\le \sum_{j=0}^{J-1}(2^j)^{-1-k_\perp}\int_{t:\; 2^jr_x-\rho<|t-s_x|<2^{j+1}r_x+\rho}|\sigma(t)|\,\dt,\]
and we apply the compact embedding $H^{1}(\Lambda)\subset L^\infty(\Lambda)$, which gives: 
\[II\le C\|\sigma\|_{H^1(\Lambda)}\sum_{j=0}^{J-1}(2^jr_x)^{-k_\perp}.\]
If $k_\perp=0$, then the summation equals $J\sim|\log(r_x)|$. On the other hand, if $k_\perp\ge 1$, the summation is bounded by $Cr_x^{-k_\perp}$. Let us consider first this second case. If $k_\perp\ge 1$ we have that: 
\begin{align*}
	\|II\|_{L^2_{\gamma,\mu}(C(\Lambda[0,\frac{L}{2}],R_0)\setminus C(\Lambda,2\rho))}^2 &\le C\|\sigma\|_{H^1(\Lambda)}\int_{C(\Lambda[0,\frac{L}{2}],R_0)\setminus C(\Lambda,2\rho))}d(x)^{2\gamma-2k_\perp}d_e(x)^{2\mu}\,\dx \\
					&\le C\|\sigma\|_{H^1(\Lambda)},
\end{align*}
where in the last step we integrated in spherical coordinates as we did for $I$. 

Finally, if $k_\perp=0$, we use that $r_x=d(x)$. Moreover, taking spherical coordinates $(r,\theta,\xi)$ in a ball containing $C(\Lambda,\frac{L}{2},R_0)$ , we have that $d(x) = r\sin(\xi)$. Hence 

\begin{align*}
	\|II\|_{L^2_{\gamma,\mu}(C(\Lambda[0,\frac{L}{2}],R_0)\setminus C(\Lambda,2\rho))}^2 &\le C\|\sigma\|_{H^1(\Lambda)}\int_{C(\Lambda[0,\frac{L}{2}],R_0)\setminus C(\Lambda,2\rho))}|\log(d(x))| d(x)^{2\gamma}d_e(x)^{2\mu}\,\dx \\  
					&\le C\|\sigma\|_{H^1(\Lambda)}^2\int_0^{2\pi}\int_0^\frac{\pi}{2}\int_{0}^{L}|\log(r\sin(\xi))| \sin(\xi)^{2\gamma+1}r^{2\gamma+2\mu+2}\,\dr\dxi\dtheta \\
					&\le C\|\sigma\|_{H^1(\Lambda)}^2\bigg\{ \int_0^{\frac{\pi}{2}}\int_{0}^{L} \log(\sin(\xi))\sin(\xi)^{2\gamma+1}r^{2\gamma+2\mu+2}\dr\dxi \\
					&\;+\int_0^{\frac{\pi}{2}}\int_{0}^{L} \log(r)\sin(\xi)^{2\gamma+1}r^{2\gamma+2\mu+2}\dr\dxi 
\bigg\} \\ 
					&\le C\|\sigma\|_{H^1(\Lambda)},
\end{align*}
where in the last step we used again that the integrals are bounded independently of $\rho$ thanks to conditions \eqref{gamma cond anisotropic} and \eqref{mu cond anisotropic}. This concludes the proof.
\end{proof}

\begin{lemma}\label{anisotropic E not so near}
Under the conditions of Proposition \ref{prop urho anisotropic not so near}, there is a constant $C$ independent of $\rho$ such that
\[\|(D^{\beta_\perp}\Gamma)*\mathcal{E}^{k_s}_\rho\|_{L^2_{\gamma,\mu}(C(\Lambda,R_0)\setminus C(\Lambda,2\rho))} \le C\|\sigma\|_{H^{k_s}(\Lambda)}.\]
\end{lemma}
\begin{proof}
	As in the term $II$ in the proof of Lemma \ref{anisotropic E near}, here it is convenient to integrate by parts, passing the derivatives with respect to $s$ from the regularized function $\phi_{1,\rho}$ to the kernel $\Gamma$. Denoting $\partial_{s_y}^{k_s-\ell-1}$ the derivative of order $k_s-\ell-1$ with respect to $s_y$, we have: 
	\begin{align*}
	\Big|\int_{\bar{\bf C}_0^\rho} D^{\beta_\perp}\Gamma(x-y)\mathcal{E}_\rho^{k_s}\,\dy\Big| &= \Big|\int_{\bar{\bf C}_0^\rho} D^{\beta_\perp}\Gamma(x-y)\sum_{\ell=0}^{k_s-1}\sigma^{(\ell)}(\rho)\phi_{1,\rho}^{(k_s-1-\ell)}(s_y-\rho)\phi_{2,\rho}(r_y)\,\dy\Big| \\
	&\le \int_{\overline{\bf C}_0^\rho} \Big|\sum_{\ell=0}^{k_s-1}\partial_{s_y}^{(k_s-1-\ell)}D_x^{\beta_\perp}\Gamma(x-y)\sigma^{(\ell)}(\rho)\phi_{1,\rho}(s_y-\rho)\phi_{2,\rho}(r_y)\,\dy\Big|, 
	\end{align*}
	where we used that $\phi_{1,\rho}\phi_{2,\rho}$ and its derivatives vanish at the boundary of ${\bf C}_0^\rho$. Now, using \eqref{estimate phi} and \eqref{embedding}, we continue with
	\begin{align*}
		&\le C\rho^{-3}\int_{\overline{\bf C}_0^\rho} \sum_{\ell=0}^{k_s-1}|x-y|^{\ell-k_s-k_\perp}\|\sigma\|_{H^{\ell+1}(\Lambda)},\dy \\
		&\le C\rho^{-3}\|\sigma\|_{H^{\ell+1}(\Lambda)}\int_{\overline{\bf C}_0^\rho} |x-y|^{-k_s-k_\perp}\,\dy 
	\end{align*}
	Moreover, since we need to take $x\in C(\Lambda,R_0)\setminus C(\Lambda,2\rho)$, we have that $|x-y|\sim d_e(x)$, using this and that $|\overline{\bf C}_0^\rho|\sim\rho^{3}$ we complete the estimate:
	\[|(D^{\beta_\perp}\Gamma)*\mathcal{E}_\rho^{k_s}| \le C\|\sigma\|_{H^{k_s}(\Lambda)} d_e(x)^{-k_s-k_\perp}.\]
	Inserting this estimate in the norm, we have
	\begin{align*}
		\|(D^{\beta_\perp}\Gamma)*&\mathcal{E}_\rho^{k_s}\|_{L^2_{\gamma,\mu}(C(\Lambda[0,\frac{L}{2}],R_0)\setminus C(\Lambda,2\rho))} \\
				&\le C\|\sigma\|_{H^{k_s}(\Lambda)}\bigg(\int_{C(\Lambda[0,\frac{L}{2}],R_0)\setminus C(\Lambda,2\rho)}d(x)^{2\gamma}d_e(x)^{2\mu-2k_s-2k_\perp}\,\dx\bigg)^\frac{1}{2},
	\end{align*}
	and as in the previous lemmas, we can complete the estimate by integrating the weights in spherical coordinates, in a ball containing $C(\Lambda[0,\frac{L}{2}],R_0)$, using conditions \eqref{gamma cond anisotropic} and \eqref{mu cond anisotropic}. 
\end{proof}

Finally, we prove the following lemma.

\begin{lemma}\label{anisotropic urho not so near ball}
Under the conditions of Proposition \ref{anisotropic urho near}, the following estimate holds: 
\[\|D^\beta u_\rho^\circ\|_{L^2_{\gamma,\mu}((B^+_{R_0}(0)\cup B^+_{R_0}(L))\setminus(B^+_{2\rho}(0)\cup B^+_{2\rho}(L))}\le C\|\sigma\|_{H^m(\Lambda)},\]
with a constant $C$ independent of $\rho$.
\end{lemma}
\begin{proof}
	The situation is quite similar to the one in Lemma \ref{anisotropic urho near ball}. Thanks to Lemma \ref{anisotropic deriv}, we have that: 
	\[\|D^\beta u_\rho^\circ\|_{L^2_{\gamma,\mu}(B^+_{R_0}(0)\setminus B^+_{2\rho}(0))} \le \|D^{\beta_\perp}(\sigma^{(k_s)})_\rho\|_{L^2_{\gamma,\mu}}(B^+_{R_0}(0)\setminus B^+_{2\rho}(0)) + \|D^{\beta_\perp}\mathcal{E}^{k_s}_\rho\|_{L^2_{\gamma,\mu}}(B^+_{R_0}(0)\setminus B^+_{2\rho}(0)).\]
	But in $B^+_{R_0}(0)$, $d(x)=d_e(x)$ so the weight is reduced to $d_e(x)^{2(\gamma+\mu)}$. 
	The estimate for the first term is almost exactly the same than the one given for Lemma \ref{anistropic sigma ks not so near}, but working only with $d_e^{2(\gamma+\mu)}$ and integrating in spherical coordinates.
	For the second term the estimate is analogous to the one for Lemma \ref{anisotropic E not so near}, with the same adaptation. 
\end{proof}

With this lemma, we have completed the proof of Theorem \ref{main theorem anisotropic}.

\section{Some extensions}\label{section extensions}

In this section we present some extensions of the results previously obtained. We discuss the main ideas that lead to these extensions, but we do not provide a detailed proof of any of them. The reader can easily fill the gaps. 

\subsection{The two dimensional case}\hfill\\
The technique applied in the previous sections can be also used for treating the two dimensional case. Naturally, cylindrical and spherical coordinates should be replaced by curvilinear cartesian and polar coordinates, respectively. 

An important issue arises, however, when trying to apply Theorem \ref{theorem SW}, since the kernel $\log(|x|)$ does not define a fractional integral. When analizing $D^\beta u_\rho$ for some $k=|\beta|>0$, this problem can be avoided by using the respresentation formula: 
\[D^\beta u_\rho^\circ(x) = \int_{C(\Lambda,\rho)} D^{\beta'}\Gamma(x-y) D^{\beta''}\sigma_\rho(y)\,\dy,\]
where $\beta = \beta'+ \beta''$ and $\beta'$ is chosen such that $|\beta'|=1$. In this case, we have that $|D^\beta u_\rho(x)| \le C|I_1(D^{\beta''}\sigma_\rho)(x)|$ and Lemma \ref{lemma Apqalpha} can be applied with $\alpha=1$ and a derivative of order $k-1$ of $\sigma_\rho$. The rest of the calculations can be carried out as in Section \ref{section isotropic}. 

However, when we want to estimate the norm of $u_\rho^\circ$ (no derivative), we need to deal with the kernel $|\log(x)|$. We can observe that $|\log(x)|\le C|x|^{\alpha-2}$ for every $\alpha<2$, which would allow us to estimate  $u_\rho^\circ$ by $I_\alpha(\sigma_\rho)$ with $\alpha$ as close to $2$ as needed. However, this is not possible. Indeed, in dimension two, the combination of conditions \eqref{Apqalphacondeta} and \eqref{Apqalphacondeta*} with \eqref{eta*} for $p<q=2$ gives the restriction: 
\[-\frac{1}{p}+\alpha<\eta<\frac{1}{2},\]
and when we take $p$ close to $2$ we are forced to take $\alpha$ close to $1$. Consequently, we obtain for $u_\rho^\circ$ the same restrictions on $\gamma$ that we have for its derivatives of first order, so we do not get the expected shift in the weight. The result that we are able to prove is the following:

\begin{thm}
	Given $\Omega\subset \R^2$ a domain of class $\mathcal{C}^{m-1,1}$, $\Lambda$ a curve strictly contained in $\Omega$ satisfying the conditions estipulated in Section \ref{section preliminaries} and $\sigma\in L^2(\Lambda)$, then, the solution $u$ to problem \eqref{problem} satisfies that $u\in L^2_\gamma(\Omega)$ and $\nabla u\in K^{m-1}_\gamma(\Omega)^2$ for every $\gamma>-\frac{1}{2}$. Moreover, there is a constant $C$ such that: 
\[\|u\|_{L^2_\gamma(\Omega)}+\|\nabla u\|_{K^{m-1}_\gamma(\Omega)^2}\le C\|\sigma\|_{L^2(\Lambda)}.\]
The result also holds for $m=2$ and $\Omega$ a convex polygon. 
\end{thm}

In a similar way we can prove the following anisotropic result: 

\begin{thm}
	Under the conditions of the previous theorem with the additional assumption that $\sigma\in H^m(\Lambda)$, we have that $u\in L^2_{\gamma,\mu}(\Omega)$ and $\nabla u\in K^{m-1}_{\gamma,\mu}(\Omega)$ for every $\gamma>-\frac{1}{2}$ and $\mu>-\frac{1}{2}$. 
\end{thm}

\subsection{Anisotropic results for curved fractures}\hfill\\
For simplicity, in Section \ref{section anisotropic} we restricted our analysis to the case where $\Lambda$ is a segment. However, in the case of a simple curve, the notions of $\beta_\perp$ and $\beta_s$ are meaningful as long as we can take cylindrical coordinates in  $C(\Lambda,R_0)$ and spherical coordinates in $B^+(0)$ and $B^+(L)$. Consequently, the Kondratiev type space $K_{\gamma,\mu}(C(\Lambda,R_0))$ is well defined.

\begin{figure}[!h]
	\includegraphics[scale=0.7]{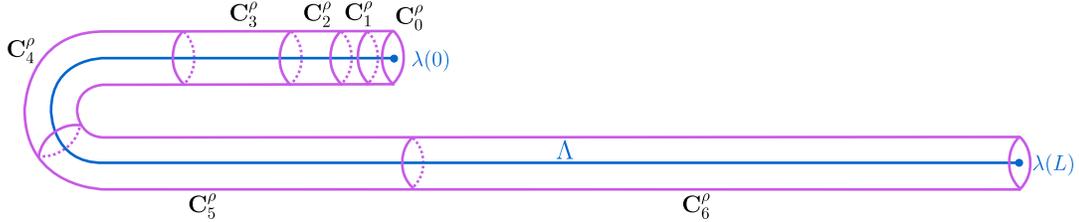}
	\caption{A general simple, open curve.}
	\label{fig curvature}
\end{figure}

The main advantage of considering a straight segment is that when dealing with the partition \eqref{decomp} it is clear that if $x\in {\bf C}_j^\rho$ for $1\le j<J$ we have that $d_e(x)\sim 2^j\rho$ and that if $x\in {\bf C}_0^\rho$ and $y\in {\bf C}_j^\rho$ for $j>1$, then $|x-y|\sim (2^j\rho)$. This might not be obvious for a general curve where a situation as the one depicted in Figure \ref{fig curvature} can occur. There, we have that there are points around the middle of $\Lambda$ (in $\mathbf{C}_5^\rho$) that are closer to the extreme point $\lambda(0)$ than some point in, for example, $\mathbf{C}_2^\rho$. However, it is easy to check that in this case the estimates for $d_e(x)$ and $|x-y|$ still hold, with proportionality constants depending on $R_0$ and $L$. 

Therefore, we conclude that Theorem \ref{main theorem anisotropic} can be extended to the general problem as it is presented in Section \ref{section preliminaries}.

\subsection{Closed curves}\hfill\\
The theory of Section \ref{section isotropic} can be applied if $\Lambda$ is a simple closed curve, smooth enough so that a set of cylindrical coordinates can be defined in a neighbouring region, e.g.: a circle. In this case, we can work with an adapted version of $\sigma_\rho$, given by: 
\[\sigma_\rho(t,\theta,s) = \phi_{2,\rho}(r)\int_{s-\rho}^{s+\rho} \sigma(t)\phi_{1,\rho}(s-t)\,\dt,\]
where $\sigma$ is extended to $\R$ periodically. 

Moreover, the anisotropic result also hold, with $\mu=0$ for every $\beta$ since there are no extreme points. 

\subsection{Polygonal fractures}\hfill\\
We now consider $\Lambda$ a polygonal with vertices $\{e_i\}_{i=0}^I$. In this case, we can write 
\[\Lambda = \bigcup_{i=1}^I \Lambda_i,\]
where $\Lambda_i$ is the segment joining the points $e_{i-1}$ and $e_i$. If $\Lambda$ is open, then $e_0$ and $e_I$ are the extreme points of $\Lambda$. If $\Lambda$ is closed, then $e_0=e_I$.

With this notation, we can define $\sigma_{i,\rho}$ the approximation of the data on $\Lambda_i$, as we defined $\sigma_\rho$ in Section \ref{section approximating problem}, and: 
\[\sigma_\rho = \sum_{i=1}^I\sigma_{i,\rho}.\]

Each $\sigma_{i,\rho}$ is supported on the cylindrical neighbourhood of $\Lambda_i$, $C(\Lambda_i,\rho)$.

Applying the results proven in the previous sections to each $\sigma_{i,\rho}$ we obtain an analogous to Theorem \ref{main theorem anisotropic}, but redefining $d_e$ as the distance to the set of vertices:  
\[d_e(x) = \min_{0\le i\le I} |x-e_i|.\]

\begin{figure}[!ht]
	\includegraphics[scale=0.8]{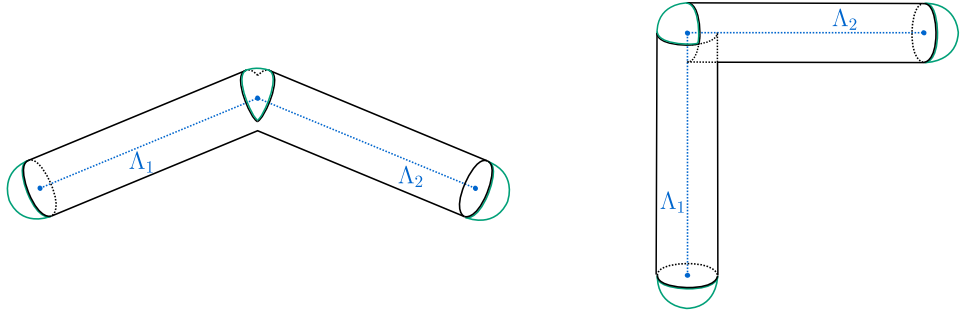}
	\caption{Examples of polygonal $\Lambda$ and the neighbourhood $B(\Lambda,\rho)$.}
	\label{figure polygonal}
\end{figure}

The only issue that needs to be addressed is that some overlapping occurs between $C(\Lambda_i,\rho)$ and $C(\Lambda_{i+1},\rho)$, as it is shown in Figure \ref{figure polygonal}. In the picture of the left the angle between $\Lambda_1$ and $\Lambda_2$ is greater that $\pi/2$. In this case, for a point $x\in C(\Lambda_1,\rho)\cap C(\Lambda_2,\rho)$ we have that $d(x,\Lambda_1)\sim d(x,\Lambda_2)$. On the contrary, in the right picture the angle between the segments is exactly $\pi/2$. This implies that there are points in $C(\Lambda_1,\rho)$ that are \emph{far} from $\Lambda_1$ but are touching $\Lambda_2$. This might seem like a problem, in particular for estimating the weighted norms of $\sigma_\rho$: $\sigma_{1,\rho}$ is defined with respect to $\Lambda_1$ but $d(x)$ can be $d(x,\Lambda_1)$ or $d(x,\Lambda_{2})$ depending on the point $x\in C(\Lambda_1,\rho)$. It is possible to check that this is actually not a problem, and the results proven in Sections \ref{section isotropic} and \ref{section anisotropic} stand, but it is quite tedious to adapt each estimate taking into account these issues.

A possible workaround is as follows. In a case like the one at the right of Figure \ref{figure polygonal} we can define $\sigma_{i,\rho}$ as: 
\[\sigma_{i,\rho}(r,\theta,s) = \phi_{2,\rho}(r) \int_{2\rho}^{L-2\rho} \sigma_i(t)\phi_{1,\rho}(s-t)\,\dt.\]
It is easy to check that Lemma \ref{lemma sigma_rho approx} still holds for this variant of $\sigma_\rho$. But now supp$(\sigma_{i,\rho})=C(\Lambda_i[\rho,L-\rho],\rho)$ and for every $x\in$supp$(\sigma_{i,\rho})$ near $e_i$ we have that $d(x,\Lambda_i)\sim d(x,\Lambda_{i+1})$ so we can estimate the weighted norms of $\sigma_{i,\rho}$ using only $d(x,\Lambda_i)$. 

Naturally, if the angle between two segments is smaller, we would need to further adapt the definition of $\sigma_{i,\rho}$. But in any case we can take
\[\sigma_{i,\rho}(r,\theta,s) = \phi_{2,\rho}(r) \int_{\kappa\rho}^{L-\kappa\rho} \sigma_i(t)\phi_{1,\rho}(s-t)\,\dt,\]
for some $\kappa$ depending on the minimal angle of the polygonal, such that for every $x\in C(\Lambda_i,\rho)\cap C(\Lambda_{i+1},\rho)$ $d(x,\Lambda_i)\sim d(x,\Lambda_{i+1})$. And with this definition, it is easy to see that Theorem \ref{main theorem anisotropic} holds.

\appendix
\section{Proof of Lemma \ref{lemma Apqalpha}}

We now prove Lemma \ref{lemma Apqalpha}. Let us begin by defining the concept of Assouad dimension.

\begin{defi}
Given $E\subset\R^n$, we denote $N_r(E)$ the smallest number of balls of radius $r$ needed for covering $E$. The Assouad dimension of $E$, denoted $\dim_A(E)$ is the infimal $\varsigma$ such that there exists a constant $C$ such that for all $0<r<R$ and $x\in E$:
	\[N_r(E\cap B(x,R))\le C\Big(\frac{R}{r}\Big)^\varsigma.\] 
\end{defi}

This definition extends naturally the behaviour of integer dimensions. If we have, for example, an $m$ dimensional manifold $E$ in $\R^n$, for some $m\in \N_0$ we can cover $E\cap B(x,R)$ with $\sim (R/r)^m$ balls of radius $r$. In particular, a dot, a smooth curve and a smooth surface have Assouad dimension $0$, $1$ and $2$ respectively. We refer the reader to \cite{Fraser} for an extensive study of the Assouad dimension. Let us just remark that the Assouad dimension is usually called \emph{the greater of all dimensions}, since it turns out to be greater than other usual dimensions. For example, the following sequence of inequalities  hold for every bounded set $E$: 
\[\dim_H(E)\le \dim_P(E)\le \overline{\dim}_B(E)\le \dim_A(E),\]
where $\dim_H$, $\dim_P$ and $\overline{\dim}_B$ are the Hausdorff, Packing and upper box dimensions, respectively. A particularly interesting example is given by the subset of the real line $E=\{0\}\cup\{\frac{1}{n}:\,n\in \N\}$, where we have: 
\[\dim_H(E)=0,\quad \overline{\dim}_B(E)=\frac{1}{2},\quad \dim_A(E)=1.\]
This simple example shows that the local nature of the Assouad dimension implies that it \emph{sees} the set as a line near the accumulation point at the origin. On the other extreme, the Hausdorff dimension, which is global, is zero for every countable set. In the middle, the box dimension captures some of the local behaviour near the origin and some of the global characteristics of the set.

For proving Lemma \ref{lemma Apqalpha} we will also need to work with Whitney decompositions, which definition we recall:

\begin{defi}[Whitney decomposition] Let $\Omega\in\R^n$ be an open set, then there exists a set of closed cubes, $\mathcal{W}=\{Q_j\}_{j\in\N}$, with edges parallel to the coordinate axis, such that  $\Omega=\cup_{j} Q_j$, satisfying:
\[\sqrt{n} \ell(Q_j)\le d(Q_j,\partial\Omega)\le 4\sqrt{n} \ell(Q_j),\]
\[\frac{1}{4}\le\frac{\ell(Q_{i})}{\ell(Q_j)}\le 4, \textrm{ if }{Q_{i}}\cap {Q_j}\neq\emptyset,\]
where $\ell(Q_j)$ is the edge length of $Q_j$.

Moreover, the cubes $\{Q_j\}$ can be assumed to be dyadic, and classified into \emph{generations}, where each generation is formed by all the cubes of a given size. Namely: the set 
\[\{Q_i^k\}_{1\le i\le W_k}=\{Q\in\mathcal{W}:\;\ell(Q)=2^{-k}\}\] 
is the $k-$th generation of cubes, with cardinal $W_k$. We also denote $W_k(B(x,R))$ the number of Whitney cubes of $k$-th generation contained in the ball $B(x,R)$. 
\end{defi}

The following result provides an estimate for $W_k(B(x,R))$. It is analogous to \cite[Lemma 6.1]{DLG_Stokes}. However, in that paper the authors consider a set contained in another set that is Alfohrs $m$-regular. Here we simplify the approach stating the lemma in terms of the Assouad dimension.

\begin{lemma}\label{L:whitneygenerations}
Let $E\subset\R^n$. Given $x\in E$, and $R>0$, there exists a constant $C$ such that
\[W_k(B(x,R))\le CR^\mu 2^{k\varsigma},\]
for every $\varsigma>\dim_A(E)$. 
\end{lemma}
\begin{proof}
	Given a cube $Q_i^k$ of $k$-th generation contained in $B(x,R)$, we define $x_i\in E$ such that $d(x_i,Q_i^k)=d(x_i,E)$. Applying the properties of a Whitney decomposition and the fact that $2^{-k}<R$, we obtain: 
	\begin{align*}
		|x-x_i| &\le d(x,Q_i^k)+\diam(Q_i^k)+d(x_i,Q_i^k) \\
		        &\le R +\sqrt{n}2^{-k}+4\sqrt{n}2^{-k} \\
						&\le C_n R,
	\end{align*}
	where $C_n>1$ depends only on the dimension $n$. Naturally, $B(x,R)\subset B(x,C_nR)$. Moreover, $B(x,C_n R)$ contains all the $x_i$ corresponding to cubes $Q_i^k$ in $B(x,R)$. 

	Applying the definition of the Assouad dimension we can cover $E\cap B(x,C_n R)$ with $N_{2^{-k}}(E\cap B(x,C_n R))$ balls $B(z_j,2^{-k})$ centered at $z_j\in E$ and with radius $2^{-k}$. Since $2^{-k}=\ell(Q_i^k)\le \frac{1}{\sqrt{n}}d(Q_i^k,E)$ there is a constant $K_n$ depending only on $n$ such that
	\[\bigcup_j B(z_j,K_n2^{-k})\]
	contains all the Whitney cubes of $k$-th generation which distance to $E$ is reached at $E\cap B(x,C_n R)$. In particular, it contains all the Whitney cubes of $k$-th generation contained in $B(x,R)$. 

	Moreover, each expanded ball $B(z_j,K_n 2^{-k})$ can pack at most $C(K_n2^{-k}/2^{-k})^n \sim C$ cubes of edge length $2^{-k}$, where $C$ is a constant depending only on $n$. Consequently the number of cubes of $k$-th generation contained in $B(x,R)$ is at most $CN_{2^{-k}}(E\cap B(x,C_n R)) \le C(C_n R/2^{-k})^\varsigma$, for  every $\varsigma>\dim_A(E)$, and the proof is finished. 
\end{proof}

Now we are finally able to prove Lemma \ref{lemma Apqalpha}. 

\begin{proof}[Proof of Lemma \ref{lemma Apqalpha}]
	Let $Q$ be a cube with edge length $\ell(Q)$ and let us denote $d_Q = d(Q,E)$. We separate the proof into two case. 
	\begin{itemize}
		\item If $\frac{\sqrt{n}}{2}\ell(Q)\le d_Q$, then for every $x\in Q$ we have that $d_Q\le d(x,E)\le 3d_Q$. Applying  this equivalence, condition \eqref{Apqalphacondalpha} and \eqref{eta*} we obtain: 
		\begin{align*}
			|Q|^{\frac{\alpha}{n}+\frac{1}{q}-\frac{1}{p}}&\bigg(\frac{1}{|Q|}\int_Q d(x,E)^{-\eta q\tau}\,\dx\bigg)^\frac{1}{q\tau}\bigg(\frac{1}{|Q|}\int_Q d(x,E)^{\eta^* p'\tau}\,\dx\bigg)^\frac{1}{p'\tau} \\
			&\le \ell(Q)^{\alpha+\frac{n}{q}-\frac{n}{p}}d_Q^{-\eta}d(x,E)^{\eta^*}\\
			&\le C d_Q^{\alpha+\frac{n}{q}-\frac{n}{p}-\eta+\eta^*} \le C
		\end{align*}
		with $C$ a constant depending only on $n$. 
		\item If $\frac{\sqrt{n}}{2}\ell(Q)> d_Q$, there is a point $x\in E$ such that the cube centered at $x$ with edges of length $3\sqrt{n}\ell(Q)$ contains $Q$. Hence, we can assume without loss of generality that $Q$'s center lies in $E$. With this assumption, we consider a Whitney decomposition of $\R^n\setminus E$ and we denote $\mathcal{W}^k_Q$ the Whitney cubes of $k$-th generation that intersects $Q$. With this notation and using that $\ell(P)\sim d(P,E)$ for every Whitney cube $P$: 
		\begin{align*}
			|Q&|^{\frac{\alpha}{n}+\frac{1}{q}-\frac{1}{p}}\bigg(\frac{1}{|Q|}\int_Q d(x,E)^{-\eta q\tau}\,\dx\bigg)^\frac{1}{q\tau}\bigg(\frac{1}{|Q|}\int_Q d(x,E)^{\eta^* p'\tau}\,\dx\bigg)^\frac{1}{p'\tau} \\
			&\le\ell_Q^{\alpha+\frac{n}{q}-\frac{n}{p}}\bigg(\frac{1}{|Q|}\sum_{k=k_0}^\infty\sum_{P\in \mathcal{W}_Q^k} \int_P d(x,E)^{-\eta q\tau}\,\dx\bigg)^\frac{1}{q\tau}\bigg(\frac{1}{|Q|}\sum_{k=k_0}^\infty\sum_{P\in \mathcal{W}_Q^k} \int_P d(x,E)^{\eta^* p'\tau}\,\dx\bigg)^\frac{1}{p'\tau} \\
			&\le\ell_Q^{\alpha+\frac{n}{q}-\frac{n}{p}}\bigg(\frac{1}{|Q|}\sum_{k=k_0}^\infty\sum_{P\in \mathcal{W}_Q^k} \ell_P^{n-\eta q\tau}\bigg)^\frac{1}{q\tau}\bigg(\frac{1}{|Q|}\sum_{k=k_0}^\infty\sum_{P\in \mathcal{W}_Q^k} \ell_P^{n+\eta^* p'\tau}\bigg)^\frac{1}{p'\tau}
		\end{align*}
	Now, we use that $\ell(P) = 2^{-k}$ for every $P\in \mathcal{W}_Q^k$. Moreover, the number of cubes in $\mathcal{W}_Q^k$ is at most the number of Whitney cubes contained in a ball centered at $x$ with radius $2\ell(Q)$. Applying this estimates we continue, taking some $\varsigma>\dim_A(E)$: 
		\begin{align*}
			&\le C\ell_Q^{\alpha+\frac{n}{q}-\frac{n}{p}}\bigg(\sum_{k=k_0}^\infty \ell(Q)^\varsigma 2^{k\varsigma} 2^{-k(n-\eta q\tau)}\bigg)^\frac{1}{q\tau}\bigg(\frac{1}{|Q|}\sum_{k=k_0}^\infty \ell(Q)^\varsigma 2^{k\varsigma} 2^{-k(n+\eta^* p'\tau)}\bigg)^\frac{1}{p'\tau} \\
			&\le C\ell_Q^{\alpha+\frac{n}{q}-\frac{n}{p}+\frac{\varsigma-n}{q\tau} + \frac{\varsigma-n}{p'\tau}}\bigg(\sum_{k=k_0}^\infty 2^{-k(-\varsigma+n-\eta q\tau)}\bigg)^\frac{1}{q\tau}\bigg(\sum_{k=k_0}^\infty  2^{-k(-\varsigma+n+\eta^* p'\tau)}\bigg)^\frac{1}{p'\tau} \\
		\end{align*}
		Thanks to restrictions \eqref{Apqalphacondeta} and \eqref{Apqalphacondeta*}, we can take $\varsigma$ close enough to $\dim_A(E)$ and $\tau$ close enough to $1$ such that the exponents of $2^{-k}$ in both summations are positive, and consequently the summations are finite. Finally, the index $k_0$ corresponds to the largest Whitney cube that intersects $Q$, so $2^{-k_0}\sim \ell(Q)$. Thus, applying \eqref{eta*} we conclude the proof:
		\begin{align*}
		&\le C\ell_Q^{\alpha+\frac{n}{q}-\frac{n}{p}+\frac{\varsigma-n}{q\tau} + \frac{\varsigma-n}{p'\tau}}2^{-k_0(\frac{-\varsigma+n}{q\tau}-\eta)}  2^{-k_0(\frac{-\varsigma+n}{p's}+\eta^* )} \\
		&\le C\ell_Q^{\alpha+\frac{n}{q}-\frac{n}{p}+\frac{\varsigma-n}{q\tau} + \frac{\varsigma-n}{p'\tau}+ \frac{-\varsigma+n}{q\tau}-\eta + \frac{-\varsigma+n}{p's}+\eta^* } \le C,
		\end{align*}
		where the constant $C$ depends only on the dimension $n$. 
	\end{itemize}
\end{proof}




\bibliographystyle{plain} \bibliography{bibdoc}
	
\end{document}